\newif\ifpictures
\picturestrue 

\documentclass[12pt]{amsart}
\usepackage{amssymb}
\usepackage{epsf}

\numberwithin{equation}{section}

\newtheorem{thm}{Theorem}
\newtheorem{prop}[thm]{Proposition}
\newtheorem{lemma}[thm]{Lemma}

\newtheorem{example}[thm]{Example}
\newtheorem{remark}[thm]{Remark}
\newtheorem{definition}[thm]{Definition}

\newenvironment{rem}{\begin{remark}\rm}{\end{remark}}
\newcounter{FNC}[page]
\def\newfootnote#1{{\addtocounter{FNC}{2}$^\fnsymbol{FNC}$%
     \let\thefootnote\relax\footnotetext{$^\fnsymbol{FNC}$#1}}}

\newcommand{\C}{\mathbb{C}}
\newcommand{\Q}{\mathbb{Q}}
\newcommand{\R}{\mathbb{R}}
\renewcommand{\P}{\mathbb{P}}

\newcommand{\calV}{\mathcal{V}}

\title[Common Transversals and Tangents in $\mathbb{P}^3$]{
  Common transversals and tangents to \\
  two lines and two quadrics in $\mathbb{P}^3$}

\author{G\'abor Megyesi}
\address{Department of Mathematics\\
        UMIST\\
        P.O. Box 88\\
        Manchester, M60 1QD \\
        England}
\email{gmegyesi@ma.umist.ac.uk}
\urladdr{http://www.ma.umist.ac.uk/gm/}

\author{Frank Sottile}
\address{Department of Mathematics\\
        University of Massachusetts\\
        Amherst, MA, 01003\\
        USA}
\email{sottile@math.umass.edu}
\urladdr{http://www.math.umass.edu/\~{}sottile}

\author{Thorsten Theobald}
\address{Zentrum Mathematik\\
        Technische Universit\"at M\"unchen\\
        D--80290 M\"unchen\\
        Germany}
\email{theobald@mathematik.tu-muenchen.de}
\urladdr{http://www-m9.mathematik.tu-muenchen.de/\~{}theobald/}

\date{03 Jun 2002}
\thanks{Research of second author supported in part by NSF grant DMS-0070494}
\begin{document}

\begin{abstract}
We solve the following geometric problem, which arises in several
three-dimensional applications in computational geometry:
For which arrangements of two lines and two spheres in ${\mathbb R}^3$
are there infinitely many lines simultaneously transversal to the two 
lines and tangent to the two spheres?

We also treat a generalization of this problem to projective quadrics:
Replacing the spheres in $\mathbb{R}^3$ by quadrics in projective space
$\mathbb{P}^3$, and fixing the lines and one general quadric, 
we give the following complete geometric description of the set of (second)
quadrics for which the 2 lines and 2 quadrics have infinitely many
transversals and tangents:
In the nine-dimensional projective space ${\mathbb P}^9$ of quadrics, 
this is a curve of degree 24 consisting of 12 plane conics,
a remarkably reducible variety.
\end{abstract}

\maketitle

\section*{Introduction}
In~\cite{theobald-lineproblems-2000}, 
one of us (Theobald) considered arrangements of $k$ lines 
and $4{-}k$ {\sl spheres} in ${\mathbb R}^3$ having infinitely many
lines simultaneously transversal to the $k$ lines and tangent to the $4{-}k$
spheres. Since for generic configurations of $k$ lines and $4{-}k$
spheres there are only finitely many common transversals/tangents,
the goal was to characterize the non-generic configurations where the discrete
and combinatorial nature of the problem is lost.
One case left open was that of two lines and two spheres. We solve 
that here.

A second purpose 
is to develop 
and present a variety of techniques from computational algebraic geometry for
tackling problems of this kind. 
Since not all our readers are familiar with these techniques, we 
explain and document these techniques, with the goal of increasing
their applicability.
For that reason, we first deal with the more general problem where we
replace the spheres in ${\mathbb R}^3$ by general quadratic surfaces
(hereafter quadrics) in complex projective 3-space $\mathbb{P}^3$.
In order to study the geometry of this problem, we fix two lines and
a quadric in general position, and describe the set of (second)
quadrics for which there are infinitely many common 
transversals/tangents in terms of an algebraic curve. 
It turns out that this set is an algebraic curve of degree~24 in the space 
$\mathbb{P}^9$ of quadrics. Factoring the ideal of this curve shows that 
it is remarkably reducible:

\begin{thm}\label{thm:24-12}
 Fix two skew lines $\ell_1$ and $\ell_2$ and a general quadric~$Q$ in
 $\mathbb{P}^3$. 
 The closure of the set of quadrics~$Q'$ for which there are infinitely many 
 lines simultaneously transversal to $\ell_1$ and $\ell_2$ and tangent to
 both $Q$ and to $Q'$ is a curve of degree $24$ in the 
 $\mathbb{P}^9$ of quadrics.
 This curve consists of $12$ plane conics.
\end{thm}

We prove this theorem by investigating the ideal defining 
the algebraic curve describing the set of (second)
quadrics. Based on this,
we prove the theorem with the aid of a computer calculation in 
the computer algebra system {\sc Singular}~\cite{SINGULAR}.
As explained in Section~\ref{se:infinitetangents}, the success
of that computation depends crucially on the preceeding analysis 
of the curve. Quite interestingly, 
there are real lines $\ell_1$ and $\ell_2$ and real quadrics
$Q$ such that all 12 components of the curve of second quadrics are real.
In general, given real lines $\ell_1$, $\ell_2$, and a real quadric $Q$, 
not all of the 12 components are defined over the real numbers.

While the beautiful and sophisticated geometry of our fundamental 
problem on lines and quadrics could be sufficient motivation to study
this geometric problem, the original motivation came from algorithmic problems
in computational geometry. 
As explained in~\cite{theobald-lineproblems-2000}, problems of this type occur in
applications where one is looking for a line or ray \emph{interacting} 
(in the sense of ``intersecting'' or 
in the sense of ``not intersecting'') with a given set of three-dimensional
bodies, if the class of admissible bodies consists of polytopes
and spheres (respectively quadrics). Concrete application classes of this
type include visibility computations with moving 
viewpoints \cite{theobald-dimacs-2001}, 
controlling a laser beam in manufacturing \cite{pellegrini-97},
or the design of \emph{envelope} data structures supporting
ray shooting queries (i.e., seeking the first sphere, if any,
met by a query ray)~\cite{aas-97}.
With regard to related treatments of the resulting algebraic-geometric 
core problems, we refer to ~\cite{MPT01,Me01,M2}. In these papers,
the question of arrangements of four (unit) spheres in $\R^3$ leading to an 
infinite number of common tangent lines is discussed from various 
viewpoints.

The present paper is structured as follows: In Section~\ref{se:plueckervisi},
we review the well-known Pl\"ucker coordinates from line geometry.
In Section~\ref{se:lines2l1q}, we characterize
the set of lines transversal to two skew lines and tangent to a quadric 
in terms of algebraic curves; we study and classify these so-called
$(2,2)$-\emph{curves}.
Then, in Section~\ref{se:infinitetangents}, we study the set of 
quadrics which (for prescribed lines $\ell_1$ and $\ell_2$) lead to
most $(2,2)$-curves. 
This includes computer-algebraic calculations, based on which we establish the
proof of Theorem~\ref{thm:24-12}.
The appendix to the paper contains annotated computer code used in the proof.
In Section~\ref{se:examples}, we give some detailed examples illustrating
the geometry described by Theorem~\ref{thm:24-12}, and complete its proof.
Finally, in Section~\ref{se:lines2l2s}, we solve the original question of 
spheres and give the complete characterization of configurations of  
two lines and two quadrics having infinitely many lines transversal to the 
lines and tangent to the quadrics. 
For a precise statement of that
characterization see Theorems~\ref{th:2l2sclassification} and~\ref{th:2l2sclass}.

\section{Pl\"ucker Coordinates\label{se:plueckervisi}}

We review the well-known \emph{Pl\"ucker coordinates} of lines in 
three-dimensional (complex)
projective space $\mathbb{P}^3$. 
For a general reference, see
\cite{hodge-pedoe-b52,clo-b96,pottmann-wallner-b2001}.  
Let $x = (x_0,x_1,x_2,x_3)^{\mathrm{T}}$ and 
$y = (y_0, y_1, y_2, y_3)^{\mathrm{T}} \in \mathbb{P}^3$ be two points spanning a
line $\ell$.
Then $\ell$ can be represented (not uniquely) by
the $4 \times 2$-matrix $L$ whose two columns are $x$ and $y$.
The Pl\"ucker vector $p=(p_{01},p_{02},p_{03},p_{12},p_{13},p_{23})^{\mathrm{T}}
\in \mathbb{P}^5$
of $\ell$ is defined by the determinants of the $2 \times 2$-submatrices
of $L$, that is, $p_{ij} := x_i y_j - x_j y_i$. 
The set $\mathbb{G}_{1,3}$ of all lines in $\mathbb{P}^3$ is called the
\emph{Grassmannian of lines in} $\mathbb{P}^3$.
The set of vectors in $\mathbb{P}^5$ satisfying the 
\emph{Pl\"ucker relation} 
\begin{equation}
  \label{eq:plueckerrelation}
  p_{01} p_{23} - p_{02} p_{13} + p_{03} p_{12} = 0
\end{equation}
is in 1-1-correspondence with $\mathbb{G}_{1,3}$. 
See, for example Theorem~11 in \S~8.6 of~\cite{clo-b96}.

A line $\ell$ intersects a line $\ell'$ in $\mathbb{P}^3$
if and only if their Pl\"ucker vectors $p$ and $p'$ satisfy
\begin{equation}
\label{eq:plueckertangenttoline}
  p_{01} p'_{23}
  - p_{02} p'_{13} 
  + p_{03} p'_{12} 
  + p_{12} p'_{03} 
  - p_{13} p'_{02} 
  + p_{23} p'_{01} = 0.
\end{equation}
Geometrically, this means that the set of lines intersecting a given line
is described by a hyperplane section of the 
Pl\"ucker quadric~(\ref{eq:plueckerrelation})
in $\mathbb{P}^5$.

In Pl\"ucker coordinates we also obtain a nice characterization
(given in \cite{M2}) of the lines tangent to a given quadric in $\mathbb{P}^3$. 
(See~\cite{theobald-lineproblems-2000} for an alternative deduction
of that characterization). 
We identify a quadric $x^{\mathrm{T}} Q x = 0$ in $\mathbb{P}^3$ with its
symmetric  $4 \times 4$-representation matrix $Q$. 
Thus the sphere with center $(c_1,c_2,c_3)^{\mathrm{T}}\in\mathbb{R}^3$ and
radius $r$ and described in $\mathbb{P}^3$ by 
$(x_1 - c_1 x_0)^2 + (x_2 - c_2 x_0)^2 + (x_3 - c_3 x_0)^2 = r^2 x_0^2$,
is identified with the matrix
\[
  \left(
    \begin{matrix}
      c_1^2 + c_2^2 + c_3^2 - r^2 & -c_1 & -c_2 & -c_3 \\
      -c_1 & 1 & 0 & 0 \\
      -c_2 & 0 & 1 & 0 \\
      -c_3 & 0 & 0 & 1 
    \end{matrix}
  \right) \, .
\]
The quadric is smooth if its representation matrix has rank 4. 
To characterize the tangent lines,
we use the second exterior power of matrices 
\begin{eqnarray*}
  \wedge^2\ \colon\ \mathbb{C}^{m\times n} & \to & 
             \mathbb{C}^{\binom{m}{2}\times\binom{n}{2}}
\end{eqnarray*}
(see~\cite[p.~145]{pottmann-wallner-b2001},\cite{M2}).
Here $\mathbb{C}^{a\times b}$ is the set of $a\times b$ matrices with complex
entries. 
The row and column indices of the resulting matrix are subsets of 
cardinality~2 of $\{1, \ldots, m\}$ and $\{1, \ldots, n\}$, respectively.
For $I=\{i_1,i_2\}$ with $1\leq i_1<i_2\leq m$ and
$J=\{j_1,j_2\}$ with $1\leq j_1<j_2\leq n$, 
 \[
   \bigl(\wedge^2 A\bigr)_{I,J}\ :=\ 
       A_{i_1,j_1}A_{i_2,j_2}-A_{i_1,j_2}A_{i_2,j_1}\,.
 \]
Let $\ell$ be a line in $\mathbb{P}^3$ and $L$ be a $4 \times 2$-matrix
representing $\ell$.
Interpreting the $6 \times 1$-matrix $\wedge^2 L$ as a vector in 
$\mathbb{P}^5$, we observe that $\wedge^2 L = p_\ell$, where $p_\ell$ is the 
Pl\"ucker vector of $\ell$.

Recall the following algebraic characterization of
tangency:
The restriction of the quadratic form to the line $\ell$ is singular, in that
either it has a double root, or it vanishes identically.
When the quadric is smooth, this implies that the line is tangent to the
quadric in the usual geometric sense.

\begin{prop} [Proposition 5.2 of~\cite{M2}]
\label{le:quadrictangentcond}
 A line $\ell \subset \mathbb{P}^3$ is tangent to a quadric $Q$ if and only
 if the Pl\"ucker vector $p_\ell$ of $\ell$ lies on the quadratic 
 hypersurface in $\mathbb{P}^5$ defined by $\wedge^2 Q$, if and only if
\begin{equation}
  \label{eq:quadrictangentcond}
  p_\ell^{\mathrm{T}} \, \bigl(\wedge^2 Q\bigr) \, p_\ell = 0.
\end{equation}
\end{prop}

For a sphere with radius $r$ and center $(c_1, c_2, c_3)^{\mathrm{T}} \in \mathbb{R}^3$ 
the quadratic form $p_\ell^{\mathrm{T}} \bigl(\wedge^2 Q\bigr) p_\ell$ is
\begin{equation}
\label{eq:plueckertangenttosphere}
  \left(
  \begin{array}{c}
    p_{01} \\ p_{02} \\ p_{03} \\ p_{12} \\ p_{13} \\ p_{23}
  \end{array}
  \right)^{\mathrm{T}}
  \left(
  \begin{array}{c@{~}c@{~}c@{~}c@{~}c@{~}c}
    c_2^2 + c_3^2 - r^2 & - c_1 c_2 & - c_1 c_3 & c_2 & c_3 & 0 \\
    - c_1 c_2 & c_1^2 + c_3^2 - r^2 & - c_2 c_3 & -c_1 & 0 & c_3 \\
    - c_1 c_3 & - c_2 c_3 & c_1^2 + c_2^2 - r^2 & 0 & -c_1 & -c_2 \\
    c_2 & -c_1 & 0 & 1 & 0 & 0 \\
    c_3 & 0 & -c_1 & 0 & 1 & 0 \\
    0 & c_3 & -c_2 & 0 & 0 & 1
  \end{array}
  \right)
  \left(
  \begin{array}{c}
    p_{01} \\ p_{02} \\ p_{03} \\ p_{12} \\ p_{13} \\ p_{23}
  \end{array}
  \right) \, .
\end{equation}

\section{Lines in ${\mathbb P}^3$ meeting 2 lines and tangent to a 
  quadric\label{se:lines2l1q}}

We work here over the ground field $\C$.
First suppose that $\ell_1$ and $\ell_2$ are lines in $\P^3$ that meet at a
point $p$ and thus span a plane $\Pi$.
Then the common transversals to $\ell_1$ and $\ell_2$ either contain $p$ or
they lie in the plane $\Pi$.
This reduces any problem involving common transversals to $\ell_1$ and
$\ell_2$ to a planar problem in $\P^2$ (or $\R^2$), and so we shall always
assume that $\ell_1$ and $\ell_2$ are skew.
Such lines have the form
 \begin{equation}\label{eq:parameter}
   \begin{array}{rcl}
     \ell_1 & = & \{ w a + x b \, : \, [w,x] \in \P^1 \} \, , \\
     \ell_2 & = & \{ y c + z d \, : \, [y,z] \in \P^1 \}\rule{0pt}{15pt}
  \end{array}
 \end{equation}
where the points $a,b,c,d\in\P^3$ are affinely independent.
We describe the set of lines meeting $\ell_1$ and $\ell_2$ that are also
tangent to a smooth quadric $Q$.
We will refer to this set as the {\it envelope} of common transversals and
tangents, or (when $\ell_1$ and $\ell_2$ are understood) simply as the
envelope of $Q$.

The parameterization of~\eqref{eq:parameter} allows us to identify each of
$\ell_1$ and $\ell_2$ with $\P^1$; the point $wa+xb\in\ell_1$ is identified
with the parameter value $[w,x]\in\P^1$, and the same for $\ell_2$.
We will use these identifications throughout this section.
In this way, any line meeting $\ell_1$ and $\ell_2$ can be identified with 
the pair $([w,x], [y,z]) \in \P^1 \times \P^1$ corresponding to its
intersections with $\ell_1$ and $\ell_2$.
By~(\ref{eq:plueckertangenttoline}), the Pl\"ucker coordinates 
$p_\ell=p_\ell(w,x,y,z)$ 
of the transversal $\ell$ passing through the points 
$w a + x b$ and $y c + z d$ are separately homogeneous of
degree~1 in each set of variables $\{w,x\}$ and $\{y,z\}$, called
\emph{bihomogeneous} of \emph{bidegree} (1,1) (see, e.g.,
\cite[\S 8.5]{clo-b96}).

By Proposition~\ref{le:quadrictangentcond},
the envelope of common transversals to $\ell_1$ and $\ell_2$
that are also tangent to $Q$ is given by the common transversals $\ell$ of
$\ell_1$  and $\ell_2$ whose Pl\"ucker coordinates $p_\ell$ additionally
satisfy $p_\ell\bigl( \wedge^2 Q \bigr) p_\ell = 0$.
This yields a homogeneous equation 
 \begin{equation}\label{eq:Fdef}
  F(w,x,y,z)\  :=\  
   p_\ell(w,x,y,z)^{\mathrm{T}} \bigl(\wedge^2 Q\bigr)p_\ell(w,x,y,z)\ =\ 0
 \end{equation}
of degree four in the variables $w,x,y,z$. 
More precisely, $F$ has the form
\begin{equation}\label{eq:22}
  F(w,x,y,z)\ =\ \sum_{i,j=0}^2 c_{ij} w^ix^{2-i}y^jz^{2-j}
\end{equation}
with coefficients $c_{ij}$, that is $F$ is bihomogeneous with bidegree $(2,2)$. 
The zero set of a (non-zero) bihomogeneous polynomial defines an algebraic
curve in $\P^1 \times \P^1$ (see the treatment of
projective elimination theory in \cite[\S 8.5]{clo-b96}).
In correspondence with its bidegree, the curve defined by $F$
is called a \emph{$(2,2)$-curve}.
The nine coefficients of this polynomial identify the set
of $(2,2)$-curves with $\mathbb{P}^8$.

It is well-known that the Cartesian product $\P^1 \times \P^1$
is isomorphic to a smooth quadric surface in
$\P^3$~\cite[Proposition~10 in~\S~8.6]{clo-b96}.
Thus the set of lines meeting $\ell_1$ and $\ell_2$ and tangent to
the quadric $Q$ is described as the intersection of two quadrics in 
a projective 3-space.
When it is smooth, this set is a genus 1 
curve~\cite[Exer.~I.7.2(d) and Exer.~II.8.4(g)]{Ha77}.
This set of lines cannot be parameterized by polynomials---only 
genus 0 curves (also called rational curves) admit such
parameterizations (see, e.g., \cite[Corollary~2 on p.268]{semple-kneebone-b59}).
This observation is the starting point for our
study of common transversals and tangents.

Let $C$ be a $(2,2)$-curve in $\mathbb{P}^1 \times \mathbb{P}^1$
defined by a bihomogeneous polynomial $F$ of bidegree~2.
The components of $C$ correspond to the irreducible factors of $F$, which
are bihomogeneous of bidegree at most $(2,2)$.
Thus any factors of $F$ must have bidegree one of $(2,2)$, $(2,1)$, $(1,1)$,
$(1,0)$, or $(0,1)$.
(Since we are working over $\C$, a homogeneous quadratic of bidegree $(2,0)$
factors into two linear factors of bidegree $(1,0)$.)
Recall (for example, \cite{clo-b96})
a point $([w_0,x_0],[y_0,z_0]) \in C\subset\mathbb{P}^1 \times \mathbb{P}^1$ is 
\emph{singular} if the gradient $\nabla F$ vanishes at that point,
$\nabla F([w_0,x_0],[y_0,z_0]) = 0$.
The curve $C$ is \emph{smooth} if it does not contain a singular point;
otherwise $C$ is \emph{singular}. 
We classify $(2,2)$-curves, up to change of coordinates on 
$\ell_1\times\ell_2$, and interchange of $\ell_1$ and $\ell_2$.
Note that an $(a,b)$-curve and a $(c,d)$-curve meet if either $ad\neq 0$ or
$bc\neq0$, and the intersection points are singular on the union of the two
curves. 

\begin{lemma}\label{le:CrClass}
 Let $C$ be a $(2,2)$-curve on $\mathbb{P}^1 \times \mathbb{P}^1$.
 Then, up to interchanging the factors of $\mathbb{P}^1\times\mathbb{P}^1$,  
 $C$ is either
 \begin{enumerate}
  \item smooth and irreducible,
  \item singular and irreducible,
  \item the union of a $(1,0)$-curve and an irreducible 
       $(1,2)$-curve,
  \item the union of two distinct irreducible $(1,1)$-curves, 
  \item a single irreducible $(1,1)$-curve, of multiplicity two, 
  \item the union of one irreducible $(1,1)$-curve, 
         one $(1,0)$-curve, and one $(0,1)$-curve,
   \item the union of two distinct $(1,0)$-curves, and two distinct
           $(0,1)$-curves, 
   \item the union of two distinct $(1,0)$-curves, and one $(0,1)$-curve of
           multiplicity two, 
   \item the union of one $(1,0)$-curve, and one $(0,1)$-curve, both of
            multiplicity two. 
 \end{enumerate}
  In particular, when $C$ is smooth it is also irreducible.
\end{lemma}

When the polynomial $F$ has repeated factors, we are in cases (5),  (8),
or~(9). 
We study the form $F$ when the quadric is reducible, that is either when $Q$
has rank 1, so that it defines a double plane, or when $Q$ has rank 2 so that
it defines the union of two planes.

\begin{lemma}\label{lem:badquadrics}
 Suppose $Q$ is a reducible quadric.
 \begin{enumerate}
  \item[(1)] If $Q$ has rank 1, then $\wedge^2Q=0$, and so the form
             $F$ in~\eqref{eq:Fdef} is identically zero.

  \item[(2)] Suppose $Q$ has rank 2, so that it defines the union of two planes 
             meeting in a line $\ell$.
             If $\ell$ is one of $\ell_1$ or $\ell_2$, then the form $F$
             in~\eqref{eq:Fdef} is identically zero.
             Otherwise the form $F$ is the square of a $(1,1)$-form, and hence
             we are in cases (5) or (9) of Lemma~\ref{le:CrClass}.
 \end{enumerate}
\end{lemma}

\begin{proof}
 The first statement is immediate.
 For the second, let $\ell'$ be a line in $\P^3$ with Pl\"ucker coordinates
 $p_{\ell'}$. 
 From the algebraic characterization of tangency of
 Proposition~\ref{le:quadrictangentcond}, 
 $p_{\ell'}^{\mathrm{T}}\bigl(\wedge^2 Q\bigr)p_{\ell'}=0$ implies that the restriction
 of the quadratic form to ${\ell'}$ either has a zero of multiplicity two, or
 it vanishes identically. 
 In either case, this implies that ${\ell'}$ meets the line $\ell$ common to
 the two planes. 
 Conversely, if ${\ell'}$ meets the line $\ell$, then
 $p_{\ell'}^{\mathrm{T}}\bigl(\wedge^2 Q\bigr)p_{\ell'}=0$.

 Thus if $\ell$ equals one of $\ell_1$ or $\ell_2$, then
 $p_{\ell'}^{\mathrm{T}}\bigl(\wedge^2 Q\bigr)p_{\ell'}=0$ for every common transversal
 $\ell'$ to $\ell_1$ and $\ell_2$, and so the form $F$ is identically zero.
 Suppose that $\ell$ is distinct from both $\ell_1$ and $\ell_2$.
 We observed earlier that the set of lines transversal to $\ell_1$ and
 $\ell_2$ that also meet $\ell$ is defined by a $(1,1)$-form $G$.
 Since the $(2,2)$-form $F$ defines the same set as does the $(1,1)$-form $G$,
 we must have that $F=G^2$, up to a constant factor.
\end{proof}

As above, let $C$ be defined by the polynomial $F$.
For a fixed point $[w,x]$, the restriction of the polynomial $F$ to
$[w,x]\times\P^1$ is a homogeneous quadratic polynomial in $y,z$.
A line passing through $[w,x]\in\ell_1$ and the point 
of $\ell_2$ corresponding to any zero of this
restriction is tangent to $Q$.
This construction gives all lines tangent to $Q$ that contain the point $[w,x]$.
We call the zeroes of this restriction the 
{\it fiber over $[w,x]$} of the projection of $C$ to $\ell_1$.

We investigate these fibers.
Consider the polynomial $F$ as a polynomial in the variables $y,z$ with
coefficients polynomials in $w,x$.
The resulting quadratic polynomial in $y,z$ has discriminant
 \begin{equation}\label{eq:disc}
   \left(\sum_{i=0}^2 c_{i1}w^ix^{2-i}\right)^2\ -\ 
   4\left(\sum_{i=0}^2 c_{i0}w^ix^{2-i}\right)
   \left(\sum_{i=0}^2 c_{i2}w^ix^{2-i}\right)\;.\smallskip
 \end{equation}

\begin{lemma} 
 If this discriminant vanishes identically, then the polynomial $F$ has
 a repeated factor.
\end{lemma}

\begin{proof}
 Let $\alpha,\beta,\gamma$ be the coefficients of $y^2,yz,z^2$ in the
 polynomial $F$, respectively.
 Then we have $\beta^2=4\alpha\gamma$, as the discriminant vanishes.
 Since the ring of polynomials in $w,x$ is a unique factorization domain, 
 either $\alpha$ differs from $\gamma$ by a constant factor, or else both
 $\alpha$ and $\gamma$ are squares.
 If $\alpha$ and $\gamma$ differ by a constant factor, then so do $\alpha$
 and $\beta$.  
 Writing $\beta=2d\alpha$ for some $d\in\mathbb{C}$, we have 
 \[
    F\ =\ \alpha y^2 + 2d\alpha yz + d^2\alpha z^2\ =\ 
          \alpha (y + d z)^2\,.
 \]
 If we have $\alpha=\delta^2$ and $\gamma=\sigma^2$ for some linear
 polynomials $\delta$ and $\sigma$, then
 \[
   F\ =\ \delta^2 y \pm 2\delta\sigma yz + \sigma^2 z^2\ =\ 
          (\delta y\pm \sigma z)^2\,.\vspace{-10pt}
 \]
\end{proof}

When $F$ does not have repeated factors, 
the discriminant does not vanish identically.
Then the fiber of $C$ over the point $[w,x]$ of $\ell_1$ consists of two
distinct points exactly when the discriminant does not vanish at $[w,x]$.
Since the discriminant has degree 4, there
are at most four fibers of $C$ consisting of a double point rather than two
distinct points.
We call the points $[w,x]$ of $\ell_1$ whose fibers consist of such double
points {\it ramification} points of the projection from $C$ to $\ell_1$.

This discussion shows how we may parameterize the curve $C$, at least locally.
Suppose that we have a point $[w,x]\in\mathbb{P}^1$ where the
discriminant~\eqref{eq:disc} does not vanish.
Then we may solve for $[y,z]$ in the polynomial
$F$ in terms of $[w,x]$.
The different branches of the square root function give
local parameterizations of the curve $C$.

\subsection{A normal form for asymmetric smooth $(2,2)$-curves\label{sec:smooth}}
Recall that for any distinct points $a_1,a_2,a_3 \in \P^1$ and any
distinct points $b_1,b_2,b_3 \in \P^1$, there exists a 
projective linear transformation (given by a regular $2 \times 2$-matrix) 
which maps $a_i$ to $b_i$, $1 \le i \le 3$~\cite{clo-b96,pottmann-wallner-b2001}.

\begin{lemma}\label{L:ramification}
 If the $(2,2)$-curve is smooth then the projection of $C$ to $\ell_1$ has four
 different ramification points.
\end{lemma}

\begin{proof}
Changing coordinates on $\ell_1$ and $\ell_2$ by a projective linear
transformation if necessary, we may 
assume that this projection to $\ell_1$ is ramified over $[w,x]=[1,0]$,
and the double root of the fiber is at $[y,z]=[1,0]$.
Restricting the polynomial $F$~\eqref{eq:22} to the fiber over
$[w,x] = [1,0]$ gives the equation
$$
  c_{22}y^2 + c_{21}yz + c_{20} z^2\ =\ 0\,.
$$
Since we assumed that this has a double root at $[y,z]=[1,0]$, we have 
$c_{21}=c_{22}=0$.

Suppose now that the projection from $C$ to $\ell_1$ is ramified at fewer than
four points.
We may assume that $[w,x]=[1,0]$ is a double root of the
discriminant~(\ref{eq:disc}), which implies that the coefficients of
$w^4$ and $w^3x$ in~(\ref{eq:disc}) vanish.
The previously derived condition $c_{21}=c_{22}=0$ implies that the
coefficient of $w^4$ vanishes and the coefficient of $w^3x$
becomes  $-4c_{20}c_{12}$.
If $c_{20}=0$, then every non-vanishing term of~(\ref{eq:22}) 
depends on $x$; hence, $x$ divides $F$, and so $C$ is reducible, and hence not
smooth.
If $c_{12}=0$ then the gradient $\nabla F$ vanishes at the point
$([1,0],[1,0])$, and so $C$ is not smooth.
\end{proof}

Suppose that $C$ is a smooth $(2,2)$-curve.
Then its projection to $\ell_1$ is ramified at four different points.
We further assume that the double points in the ramified fibers project to
at least 3 distinct points in $\ell_2$.
We call such a smooth $(2,2)$-curve {\it asymmetric}.
The choice of this terminology will become clear in Section~\ref{se:examples}.
We will give a normal form for such asymmetric smooth curves.

Hence, we may assume that three of the ramification points are 
$[w,x]=[0,1]$, $[1,0]$, and $[1,1]$, and the double points in these
ramification fibers occur at $[y,z]=[0,1]$, $[1,0]$, and $[1,1]$,
respectively.
As in the proof of Lemma~\ref{L:ramification}, 
the double point at $[y,z]=[1,0]$ in the fiber over $[w,x]=[1,0]$
implies that $c_{21}=c_{22}=0$.
Similarly, the double point at $[y,z]=[0,1]$ in the fiber over $[w,x]=[0,1]$
implies that $c_{00}=c_{01}=0$.
Thus the polynomial $F$~\eqref{eq:22} becomes
$$
  c_{20}w^2z^2+c_{10}wxz^2
              +c_{11}wxyz 
              +c_{12}wxy^2+c_{02}x^2y^2
$$
Restricting $F$ to the fiber of $[w,x] = [1,1]$ gives
$$
  c_{10}z^2 + c_{20}z^2\ +\ c_{11}yz\ +\ c_{02}y^2 + c_{12}y^2\,.
$$
Since this has a double root at $[y,z]=[1,1]$, we must have
$$
  -\frac{1}{2} c_{11}\ =\ c_{10} + c_{20}\ =\ c_{02} + c_{12}\,.
$$
Dehomogenizing (setting $c_{11}=-2$) and letting $c_{20}:=s$ and $c_{02}:=t$
for some $s,t \in \mathbb{C}$, we obtain the following theorem.

\begin{thm}\label{thm:gen22}
  After projective linear transformations in $\ell_1$ and $\ell_2$, an
  asymmetric smooth $(2,2)$-curve is the zero set of a polynomial
 \begin{equation}\label{eq:gen22}
  s w^2z^2 \,+\, (1{-}s) wxz^2 \,-\, 2 wxyz \,+\, 
  (1{-}t) wxy^2 \,+\, t x^2 y^2\,,
 \end{equation}
 for some $(s,t)\in\mathbb{C}^2$ satisfying
 \begin{equation}\label{eq:discrimin}
   s t (s{-}1) (t{-}1) (s{-}t) \ \neq\ 0\,. 
 \end{equation}
\end{thm}

We complete the proof of Theorem~\ref{thm:gen22}.
The discriminant~(\ref{eq:disc}) of the polynomial~(\ref{eq:gen22}) is
 \[
   4 w x (w {-} x) \left( s(t{-}1) w \,-\, t(s{-}1) x\right)\,,
 \]
which has roots at $[w,x]=[0,1], [1,0], [1,1]$, and 
$\alpha=[t(s-1),\, s(t-1)]$.
Since we assumed that these are distinct, the fourth point $\alpha$ must
differ from the first three, which implies that $(s,t)$
satisfies~\eqref{eq:discrimin}. 
The double point in the fiber over $\alpha$ occurs at 
$[y,z]=[s-1,t-1]$. 
This equals a double point in another ramification fiber only for values
of the parameters not allowed by~(\ref{eq:discrimin}). 

\begin{rem}
 These calculations show that smooth $(2,2)$-curves exhibit the following
 dichotomy.
 Either the double points in the ramification fibers project to four distinct
 points in $\ell_2$ or to two distinct points.
 They must project to at least two points, as there are at most two points in
 each fiber of the projection to $\ell_2$.
 We showed that if they project to at least three, then they project to four.
\end{rem}

We compute the parameters $s$ and $t$ from the intrinsic geometry of the curve
$C$.  
Recall the following definition of the cross ratio (see, for
example~\cite[\S 1.1.4]{pottmann-wallner-b2001}).

\begin{definition}\label{de:crossratio}
For four points $a_1, \ldots, a_4 \in \P^1$ with $a_i = [\alpha_i, \beta_i]$,
the \emph{cross ratio} of $a_1, \ldots, a_4$ is the point of\/ $\P^1$ 
defined by
\[
  \left[ 
    \frac{ \det \left( \begin{matrix}
      \alpha_1 & \alpha_4 \\
      \beta_1 & \beta_4
    \end{matrix} \right)}
    { \det \left( \begin{matrix}
      \alpha_1 & \alpha_3 \\
      \beta_1 & \beta_3
    \end{matrix} \right)} \, , \,
    \frac{ \det \left( \begin{matrix}
      \alpha_2 & \alpha_4 \\
      \beta_2 & \beta_4
    \end{matrix} \right)}
    { \det \left( \begin{matrix}
      \alpha_2 & \alpha_3 \\
      \beta_2 & \beta_3
    \end{matrix} \right)} \right]\  .
\]
If the points are of the form $a_i = [1,\beta_i]$, this simplifies
to 
\[
    \left[\frac{\beta_4 - \beta_1}{\beta_3 - \beta_1} \, , \,
    \frac{\beta_4 - \beta_2}{\beta_3 - \beta_2}\right] \, .
\]
The cross ratio of four points  $a_1,a_2,a_3,a_4 \in \P^1$ remains
invariant under any projective linear transformation.
\end{definition}

The projection of $C$ to $\ell_1$ is ramified over the points 
$[w,x]=[0,1], [1,0], [1,1]$ and $\alpha=[t(s-1),\, s(t-1)]$.
The cross ratio of these four (ordered) ramification points
is $[t(s{-}1),\, s(t{-}1)]$.
Similarly, the cross ratio of the four (ordered) double points in the
ramification fibers  is $[s{-}1,\, t{-}1]$. 

This computation of cross ratios allows us to compute the normal form of an
asymmetric smooth $(2,2)$-curve.
Namely, let $a_1,a_2,a_3$, and $a_4$ be the four ramification points of the
projection of $C$ to $\ell_1$ and $b_1,b_2,b_3$, and $b_4$ be the images in
$\ell_2$ of the corresponding double points.
Let $\gamma_1$ be the cross ratio of the four points $a_1,a_2,a_3$, and $a_4$
(this is well-defined, as cross ratios are invariant under projective linear
transformation).
Similarly, let $\gamma_2$ be the cross ratio of the points $b_1,b_2,b_3$,
and $b_4$. 
For four distinct points, the cross ratio  is an element of
$\mathbb{C} \setminus \{0,1\}$, so we express $\gamma_1,\gamma_2$ as complex
numbers. 
The invariance of the cross ratios yields the conditions on $s$ and $t$
\[
  \frac{s (t{-}1)}{t(s{-}1)}\ =\ \gamma_1 \quad \text{ and } \quad
  \frac{t{-}1}{s{-}1}\ =\ \gamma_2 \, .
\]
Again, since $\gamma_1,\gamma_2\in \mathbb{C} \setminus \{0,1\}$, these two
equations have the unique solution
\[
  s\ =\ \frac{\gamma_1 (\gamma_2 - 1)}
      {\gamma_2 (\gamma_1 - 1)} \quad \text{ and } \quad
  t\ =\ \frac{\gamma_2 - 1}{\gamma_1 - 1} \, .
\]

\section{Proof of Theorem~1}\label{se:infinitetangents}

We characterize the quadrics $Q$ which generate the same envelope of
tangents as a given quadric. 
A symmetric $4 \times 4$ matrix has 10 independent entries which identifies
the space of quadrics with $\P^9$.
Central to our analysis is a map $\varphi$ defined for almost all quadrics
$Q$. 
For a quadric $Q$ (considered as a point in $\P^9$) whose associated
$(2,2)$-form~\eqref{eq:Fdef} is \emph{not} identically zero, we let
$\varphi(Q)$ be this $(2,2)$-form, considered as a point in $\P^8$.
With this definition, we see that the Theorem~\ref{thm:24-12} is concerned
with the fiber $\varphi^{-1}(C)$, where $C$ is the $(2,2)$-curve associated to
a general quadric $Q$.
Since the domain of $\varphi$ is 9-dimensional while its range is
8-dimensional, we expect each fiber to be 1-dimensional.

We will show that every smooth $(2,2)$ curve arises as $\varphi(Q)$ for some
quadric $Q$.
It is these quadrics that we meant by general quadrics in the statement of
Theorem~\ref{thm:24-12}. 
This implies that Theorem~\ref{thm:24-12} is a consequence of the following
theorem.\medskip 
 
\begin{thm}~\label{thm:12conics}
 Let $C\in\P^8$ be a smooth $(2,2)$-curve.
 Then the closure $\overline{\varphi^{-1}(C)}$ in $\P^9$ of the fiber  of
 $\varphi$ is a curve of degree $24$ that is the union of\/ $12$ plane conics.  
\end{thm}

We prove Theorem~\ref{thm:12conics} by computing the ideal $J$ of
the fiber $\varphi^{-1}(C)$.
Then we factor $J$ into several ideals, which corresponds
to decomposing the curve of degree~24 into the union of several curves.
Finally, we analyze
the output of these computations by hand to prove the desired result.

Our initial formulation of the problem gives an ideal $I$ that not only
defines the fiber of $\varphi$, but also the subset of $\P^9$ where
$\varphi$ is not defined. 
We identify and remove this subset from $I$ 
in several costly auxiliary computations that are performed in the 
computer algebra system {\sc Singular}~\cite{SINGULAR}. 
It is only after removing the excess components that we obtain the ideal $J$ 
of the fiber $\varphi^{-1}(C)$.

Since we want to analyze this decomposition for \emph{every} 
smooth $(2,2)$-curve, we must treat the representation of $C$ as symbolic
parameters. 
This leads to additional difficulties, which we circumvent.
It is quite remarkable that the computer-algebraic calculation
succeeds and that it is still possible to analyze its result. 

In the following, we assume that $\ell_1$ is the 
$x$-axis. 
Furthermore, we may apply a projective linear transformation and 
assume without loss of generality that $\ell_2$ is the $yz$-line 
at infinity.
Thus we have 
\begin{eqnarray*}
  \ell_1 &=& \{ (w,x,0,0)^{\mathrm{T}} \in \P^3 \, \colon \, [w,x] \in \P^1 \} \, , \\
  \ell_2 &=& \{ (0,0,y,z)^{\mathrm{T}} \in \P^3 \, \colon \, [y,z] \in \P^1 \} \, .
\end{eqnarray*}
Hence, in Pl\"ucker coordinates, the lines intersecting $\ell_1$ and $\ell_2$
are given by
\begin{equation}
  \label{eq:transversalsl1l2}
  \{ (0,wy,wz,xy,xz,0)^{\mathrm{T}} \in \P^5 \, : \, [w,x], [y,z] \in \P^1 \} \, .
\end{equation}
By Proposition~\ref{le:quadrictangentcond},
the envelope of common transversals to $\ell_1$ and $\ell_2$
that are also tangent to $Q$ is given by those lines
in~(\ref{eq:transversalsl1l2}) which additionally satisfy
\begin{equation}
\label{eq:tangentQ}
  (0,wy,wz,xy,xz,0) \: \bigl(\wedge^2\!Q\bigr) \; (0,wy,wz,xy,xz,0)^{\mathrm{T}} = 0 \, .
\end{equation}

A quadric $Q$ in $\mathbb{P}^3$ is given by the quadratic form associated to
a symmetric $4 \times 4$-matrix
 \begin{equation}\label{eq:symmat}
  Q\ :=\ \left(
  \begin{matrix}
    a & b & c & d \\
    b & e & f & g \\
    c & f & h & k \\
    d & g & k & l
  \end{matrix}
  \right)\ .
 \end{equation}
In a straightforward approach the ideal $I$ of quadrics giving a 
general $(2,2)$-curve $C$ is obtained by first 
expanding the left hand side of~(\ref{eq:tangentQ}) into
 \begin{equation}\label{eq:G}
  \begin{array}{c}  
    (el{-}g^2)x^2z^2\ +\ 2(bl{-}dg)wxz^2\ +\ (al{-}d^2)w^2z^2\\
   +\ 2(ek{-}gf)x^2yz\ +\ 2(2bk{-}cg{-}df)wxyz\ +\ 2(ak{-}dc)w^2yz
    \rule{0pt}{14pt}\\
   +\ (eh{-}f^2)x^2y^2\ +\ 2(bh{-}cf)wxy^2\ +\ (ah{-}c^2)w^2y^2\,.
     \rule{0pt}{14pt}
  \end{array}
 \end{equation}
We equate this $(2,2)$-form with the general 
$(2,2)$-form~(\ref{eq:22}), as points in $\P^8$.
This is accomplished by requiring that they are proportional,
or rather that the $2\times 9$ matrix of their coefficients
\[
  \left( \begin{matrix}
    c_{00} & c_{10} & c_{20} & 
    c_{01} & c_{11} & \ldots & c_{22} \\
    el - g^2 & 2(bl - dg) & al-d^2 &
    2(ek - gf) & 2 (2bk-cg-df) & \ldots & ah-c^2
    \end{matrix} \right) 
\] 
has rank 1.
Thus the ideal $I$ is generated by the $\binom{9}{2}$ minors of
this coefficient matrix.

With this formulation, the ideal $I$ will define the fiber $\varphi^{-1}(C)$
as well as additional, excess components that we wish to exclude.
For example, the variety in $\P^9$ defined by the vanishing of the
entries in the second row of this matrix will lie in the variety 
$I$, but these points are not those that we seek.
Geometrically, these excess components are precisely where the map
$\varphi$ is not defined.
By Lemma~\ref{lem:badquadrics}, we can identify three of these excess
components, those points of $\P^9$ corresponding to rank 1 quadrics,
and those corresponding to rank 2 quadrics consisting of the union of two planes
meeting in either $\ell_1$ or in $\ell_2$.
The rank one quadrics have ideal $E_1$ generated by the entries of the matrix
$\wedge^2Q$, the rank 2 quadrics whose planes meet in $\ell_1$ have ideal
$E_2$ generated by $a,b,c,d,e,f,g$, and those whose plane meets in $\ell_2$
have ideal $E_3$ generated by $c,d,f,g,h,k,l$.

We remove these excess components from our
ideal $I$ to obtain an ideal $J$ whose set of zeroes contain the fiber
$\varphi^{-1}(C)$. 
After factoring $J$ into its irreducible components, we will observe that
$\varphi$ does not vanish identically on any component of $J$, completing the
proof that $J$ is the ideal of $\varphi^{-1}(C)$, and also the proof of
Theorem~\ref{thm:12conics}. 
\medskip

Since $c_{00}, c_{10}, \ldots, c_{22}$ have to be treated as parameters,
the computation should be carried out over 
the function field $\mathbb{Q}(c_{00},c_{10}, \ldots, c_{22})$.
That computation is infeasible.
Even the initial computation of a Gr\"obner basis for the ideal $I$ (a
necessary prerequisite) did not terminate in two days.
In contrast, the computation we finally describe terminates in 7 minutes on
the same computer.
This is because the original computation in 
$\mathbb{Q}(c_{00},c_{10},\ldots,c_{22})[a,b,\dotsc,l]$ involved too many
parameters. 
\medskip

We instead use the 2-parameter normal form~\eqref{eq:gen22} for
asymmetric smooth $(2,2)$-curves.
This will prove Theorem~\ref{thm:12conics} in the case when $C$ is an
asymmetric smooth $(2,2)$-curve.
We treat the remaining cases of symmetric smooth $(2,2)$-curves in
Section~\ref{se:examples}. 
As described in Section~\ref{sec:smooth}, by changing the coordinates on
$\ell_1$ and $\ell_2$, every asymmetric smooth $(2,2)$-curve can be
transformed into one 
defined by a polynomial in the family~\eqref{eq:gen22}. 
Equating the $(2,2)$-form~\eqref{eq:G} with the 
form~\eqref{eq:gen22} gives the ideal $I$ generated by the 
following polynomials:
 \begin{equation}\label{eq:quad-gens}
   el-g^2  \,,\ 
   ek-gf \,,\ 
   ak-dc \,,\ 
   ah-c^2  \,,\
 \end{equation}
and the ten $2\times 2$ minors of the coefficient matrix:
 \begin{equation}\label{eq:coeff-matrix}
  M\ :=\   \left(\begin{array}{ccccc}
             s    &    1-s   &      -2      &    1-t   &   t   \\   
           al-d^2 & 2(bl-dg) & 2(2bk-cg-df) & 2(bh-cf) & eh-f^2 
      \end{array}\right)\ .
 \end{equation}

This ideal $I$ defines the same three excess components as before, and we must
remove them to obtain the desired ideal $J$.
Although the ideal $I$ should be treated in
the ring $S := \mathbb{Q}(s,t)[a,b,c,d,e,f,g,h,k,l]$, 
the necessary calculations are infeasible even in this ring,
and we instead work in 
subring $R := \mathbb{Q}[a,b,c,d,e,f,g,h,k,l][s,t]$.
In the ring $R$, the ideal $I$ is homogeneous in the set of 
variables $a,b,\ldots,l$, thus
defining a subvariety of $\mathbb{P}^9\times\mathbb{C}^2$.
The ideals $E_1$, $E_2$, and 
$E_3$ describing the excess components
satisfy $E_j \supset I$, $1 \le j \le 3$.

A {\sc Singular} computation shows that $I$ is a five-dimensional subvariety
of $\mathbb{P}^9\times\mathbb{C}^2$
(see the Appendix for details). 
Moreover, the dimensions of the three excess components 
are 5, 4, and 4, respectively. In fact, it is quite easy to see
that $\dim~E_2=\dim~E_3=4$ as both ideals are defined by 7 independent linear
equations. 
\medskip

We are faced with a geometric situation of the following form.
We have an ideal $I$ whose variety contains an excess component defined by an
ideal $E$ and we want to compute the ideal of the difference
 \[
    \calV(I) - \calV(E)\,,
 \]
here, $\calV(K)$ is the variety of an ideal $K$.
Computational algebraic geometry gives us an effective method to accomplish
this, namely saturation.
The elementary notion is that of the ideal quotient $(I : E)$, which is
defined by
 \[
    (I:E)\ :=\ \{f\in R \mid fg\in I \mbox{ for all } g\in E\}\,.
 \]
Then the saturation of $I$ with respect to $E$ is
 \[
    (I:E^\infty)\ :=\ \bigcup_{n=1}^{\infty} (I : E^n) \,.
 \]
The least number $n$ such that $(I:E^\infty)=(I:E^n)$ is called the
saturation exponent.

\begin{prop}[{\cite[\S 4.4]{clo-b96}} or {\cite[\S 15.10]{eisenbud-b95}}
 or the reference manual for {\sc Singular}]
 Over an algebraically closed field, 
\[    \calV(I:E^\infty)\ =\ 
         \overline{\calV(I) - \calV(E)}\,.   \]
\end{prop}

A {\sc Singular} computation shows that the saturation exponent of the first
excess ideal $E_1$ in $I$ is 1, and so the ideal quotient suffices to
remove the excess component $\calV(E_1)$ from $\calV(I)$.
Set $I':=(I:E_1)$, an ideal of dimension 4.
The excess ideals $E_2$ and $E_3$ each have saturation exponent 4 in $I_1$, and
so we saturate $I'$ with respect to each to obtain an ideal
$J:=(\, (I':E_2^\infty) : E_3^\infty)$, which has dimension 3 in
$\P^9\times\C^2$.

To study the components of $\calV(J)$, we first apply the factorization
Gr\"obner basis algorithm to $J$,
as implemented in the {\sc Singular} command
{\tt facstd}
(see~\cite{mmn-89} or the reference manual of {\sc Singular}).
This algorithm takes two arguments, an ideal $I$ and a list
$L=f_1,\ldots,f_n$ of polynomials.
It proceeds as in the usual Buchberger algorithm to compute a Gr\"obner basis
for $I$, except that whenever it computes a Gr\"obner basis element $G$ that
it can factor, it splits the calculation into subcalculations, one for each
factor of $G$ that is not in the list $L$, adding that factor to the 
Gr\"obner basis for the corresponding subcalculation.
The output of {\tt facstd} is a list $I_1,I_2,\dotsc,I_m$ of ideals with the
property that
 \[
   \overline{\bigcup_{j=1}^m\calV(I_j)\ -\calV(f_1\dotsb f_n)}\ =\ 
   \overline{\calV(I)-\calV(f_1\dotsb f_n)}\ .
 \]
Thus, the zero set of $I$ coincides with the union of zero sets of the
factors $I_j$, in the region where none of the polynomials in the list $L$
vanish. 
In terms of saturation, this is
 \begin{equation}  
  \mathrm{rad}(I_1 \dotsb I_m  : (f_1f_2\dotsb f_n)^{\infty})\ =\ 
  \mathrm{rad}(I : (f_1f_2\dotsb f_n)^{\infty})
 \end{equation}
where $\mathrm{rad}(K)$ denotes the radical of an ideal $K$.
Some of the ideals $I_j$ may be spurious in that $\calV(I_j)$
is already contained in the union of the other $\calV(I_i)$.
\smallskip

We run {\tt facstd} on the ideal $J$ with the list of polynomials $s$, $t$,
$s{-}1$, $t{-}1$, and $s{-}t$, and obtain seven components
$J_0,J_1,\ldots,J_6$.
The components $J_1,\dotsc,J_6$ each have dimension 3, while 
the component $J_0$ has dimension 2.
Since $\calV(J_0)$ is contained in the union
of the $\calV(J_1),\dotsc,\calV(J_6)$, 
it is spurious and so we disregard it.

We now, finally, change from the base ring $R$ to the base ring $S$, and
compute with the parameters $s,t$.
There, $J$ defines an ideal of dimension 1 and degree 24 in the 9-dimensional
projective space over the field $\Q(s,t)$.
As we remarked before, we have that $\calV(J)\supset\varphi^{-1}(C)$.
The factorization of $J$ into $J_1,\dotsc,J_6$ remains valid over $S$.
The reason we did not compute the factorization over $S$ is that {\tt facstd}
and the saturations were infeasible over $S$, and the standard arguments from
computational algebraic geometry we
have given show that it suffices to compute without parameters, as long as
care is taken when interpreting the output.

Each of the factors $J_i$ has dimension 1 and degree 4.
Moreover, each ideal contains a homogeneous quadratic polynomial in the
variables $k,l$ which must factor over some field extension of $\Q(s,t)$.
In fact, these six quadratic polynomials all factor over the field
$\mathbb{Q}(\sqrt{s},\sqrt{t})$.
For example, two of the $J_i$ contain the polynomial
$(s-1)k^2-2kl-l^2$, which is the product
 \[
  \left( (\sqrt{s}{+}1)k \,+\, l\right) 
  \left((\sqrt{s}{-}1)k \,-\, l\right) \, .
 \]
For each ideal $J_i$, the factorization of the quadratic polynomial
induces a factorization of $J_i$ into two ideals
$J_{i1}$ and $J_{i2}$.
Inspecting a Gr\"obner basis for each ideal shows that each defines a plane
conic in $\P^9$.
Thus, over the field $\Q(\sqrt{s},\sqrt{t})$, $J$ defines 12 plane
conics. 

Theorem~\ref{thm:12conics} is a consequence of the following two observations.
\begin{enumerate}
\item[(1)] The factorization of $J$ gives 12 distinct components for all
          values of the parameters $s,t$ satisfying~\eqref{eq:discrimin}.
\item[(2)] The map $\varphi$ does not vanish identically on any of the
           components $\calV(J_{ij})$ for values of the parameters $s,t$
           satisfying~\eqref{eq:discrimin}. 
\end{enumerate}
By (1), no component of $J$ is empty for any $s,t$
satisfying~\eqref{eq:discrimin} and thus, for every asymmetric $(2,2)$-curve
$C$, there is a quadric $Q$ with $\varphi(Q)=C$.
Also by (1), $J$ has exactly 12 components with each a plane conic, for any $s,t$
satisfying~\eqref{eq:discrimin}, and by (2),
$\calV(J)=\overline{\varphi^{-1}(C)}$.


\section{Symmetric smooth $(2,2)$-curves}\label{se:examples}

We investigate smooth 
curves $C$ whose double points in the ramified fibers over $\ell_1$ have 
{\it only} two distinct projections to $\ell_2$.
Assume that the ramification is at the points $[w,x]=[1,1], [1,-1], [1,s]$,
and at $[1,-s]$, for some $s\in\mathbb{C} \setminus \{0,\pm1\}$ with the double points
in the fibers at  $[y,z]=[1,0]$ for the first two and at $[0,1]$ for the
second two. 
Since the points $[1,1], [1,-1], [1,s]$, and $[1,-s]$ have cross ratio
\[
   \left[ \frac{1+s}{1-s},\ \frac{1-s}{1+s}\right]\ =\ 
    \left[  1,\ \frac{(1-s)^2}{(1+s)^2}\right]\,,
\]
we see that all cross ratios in $\mathbb{P}^1 \setminus \{[1,0],[0,1],[1,1]\}$
are obtained for some $s\in\mathbb{C} \setminus \{0,\pm1\}$.
Thus our choice of ramification results in no loss of generality. 

As in Section~\ref{se:lines2l1q}, these conditions give equations on the
coefficients $c_{ij}$ of the general $(2,2)$-curve~\eqref{eq:22}:
\[
  \begin{array}{rrr}
    c_{00}+c_{10}+c_{20}\ =\ 0,&c_{01}+c_{11}+c_{21}\ =\ 0,&
    c_{00}-c_{10}+c_{20}\ =\ 0,\\
    c_{01}-c_{11}+c_{21}\ =\ 0,& c_{02}+c_{12}s+c_{22}s^2\ =\ 0,&
    c_{01}+c_{11}s+c_{21}s^2\ =\ 0,\\ 
    c_{02}-c_{12}s+c_{22}s^2\ =\ 0,& c_{01}-c_{11}s+c_{21}s^2\ =\ 0.
  \end{array}
\]
These equations have the following consequences
\[
   0\ =\ c_{21}\ =\ c_{01}\ =\ c_{12}\ =\ c_{11}\ =\ c_{10}\ 
       =\ c_{02} + c_{22}s^2\ =\ c_{00} + c_{20}\ .
\]
Hence after normalizing by setting $c_{20}=1$, the $(2,2)$-form~\eqref{eq:22}
becomes 
\[
  (x^2-w^2)y^2 + c_{22}(x^2-s^2w^2)z^2\,.
\]
While the choice of ramification points $[1,1], [1,-1], [1,s], [1,-s]$ fixes
the parameterization of $\ell_1$, the double points in the fibers of
$[1,0]$ and $[0,1]$ do not fix the parameterization of $\ell_2$.
Thus we are still free to scale the $z$-coordinate.
We normalize this equation setting $c_{22}=\pm 1$.
We do not simply set $c_{22}=1$ because that misses an important
real form of the polynomial. 
This normalization gives
 \begin{equation}\label{eq:spec22}
  (x^2-w^2)y^2\ \pm\ (x^2-s^2w^2)z^2\ \ =\ \ 
   (y^2 \pm z^2)x^2  - (y^2 \pm s^2z^2)w^2 \,.
 \end{equation}
This shows the equation to be symmetric under the
involution $[w,x]\leftrightarrow [\sqrt{\mp1}z,y]$.
This symmetry is the source of our terminology for the two classes of 
$(2,2)$-curves.
Also, if $s \not\in \{ \pm 1, 0\}$, then this is the equation of a
smooth $(2,2)$-curve.
With the choice of sign $(-)$, which we call the curve $C(s)$.

Note that~\eqref{eq:spec22}  is real if $s$ either is real or is purely imaginary
($s\in\R\sqrt{-1}$ ).
We complete the proof of Theorem~\ref{thm:24-12} with the following result for
symmetric $(2,2)$-curves.
\bigskip

\begin{thm}\label{thm:realfamily}
 For each $s\in\C \setminus \{\pm 1, 0\}$, the closure of the fiber 
 $\varphi^{-1}(C(s))$
 consists of 12 distinct plane conics.
 When $s\in\R$ or $s\in\R\sqrt{-1}$ and we use the real form
 of~\eqref{eq:spec22} with the plus sign $(+)$, then exactly 4 of these 12
 components will be real. 
 If we use the real form of~\eqref{eq:spec22} with the minus sign $(-)$,
 then if $s\in\R$, all 12 components will be real, but if $s\in\R\sqrt{-1}$,
 then exactly 4 of these 12 components will be real. 
\end{thm}

\begin{proof} Our proof follows the proof of Theorem~\ref{thm:12conics} almost exactly, but
 with significant simplifications and a case analysis.
 Unlike the proof described in Section~\ref{se:infinitetangents}, we do not give
 annotated {\sc Singular} code in an appendix, but rather supply such 
 annotated {\sc Singular} code on the web 
 page\newfootnote{{\tt http://www.math.umass.edu/\~{}sottile/pages/2l2s/}}.

 The outline is as before, except that we work over the ring of parameters
 $\Q(s)$, and find no extraneous components when we factor the ideal into
 components. 
 We formulate this as a system of equations, remove the same three excess
 components, and then factor the resulting ideal.
 We do this calculation four times, once for each choice of sign $(\pm)$
 in~\eqref{eq:spec22}, and for $s\in\R$ and $s\in\R\sqrt{-1}$. 
 Examining the output proves the result.
\end{proof}

We consider in some detail four cases of the geometry studied in
Section~\ref{se:lines2l1q}, which correspond to the four real cases of
Theorem~\ref{thm:realfamily}.
As in Section~\ref{se:lines2l1q}, let $\ell_1$ be the $x$-axis and $\ell_2$ be
the $yz$-line  at infinity. 
Viewed in $\mathbb{R}^3$, lines transversal to $\ell_1$ and $\ell_2$ are the
set of lines perpendicular to the $x$-axis.
For a transversal line $\ell$, the coordinates $[y,z]$ of the point 
$\ell \cap \ell_2$ can be interpreted as the slope of $\ell$ 
in the two-dimensional plane orthogonal to the $x$-axis.

Consider real quadrics given by an equation of the form
 \begin{equation}\label{eq:real-quad}
   x^2+(y-y_0)^2\pm z^2\ =\ 1\,.
 \end{equation}
The quadrics with the plus $(+)$ sign are spheres with center
$(0,y_0,0)^{\mathrm{T}}$ 
and radius $1$, and those with the minus $(-)$
sign are hyperboloids of one sheet.
When $|y_0|>1$ the quadric does not meet the $x$-axis.
We look at four families of such quadrics: spheres and hyperboloids that meet
and do not meet the $x$-axis.
We remark that
quadrics which are tangent to the $x$-axis give singular $(2,2)$-curves.

First, consider the resulting $(2,2)$-curve
$$
  (x^2\, -\, w^2)y^2\, \pm\, (x^2\,-\,(1\,-\,y_0^2)w^2)z^2\,.
$$
Thus we see that these correspond to the case $s=\sqrt{1-y_0^2}$ in the
parameterization of symmetric $(2,2)$-curves given above~\eqref{eq:spec22},
while in~\eqref{eq:real-quad} and~\eqref{eq:spec22} the signs $(\pm)$
correspond.

Figures~\ref{fig:one} and~\ref{fig:two} display pictures of these 
four quadrics, together with the $x$-axis, some tangents perpendicular
to the $x$-axis, and the curve on the quadric where the lines are tangent. 
\ifpictures
\begin{figure}[htb]
$$
  \begin{picture}(150,135)(0,-10)
   \put(0,10){\epsfxsize=2in\epsffile{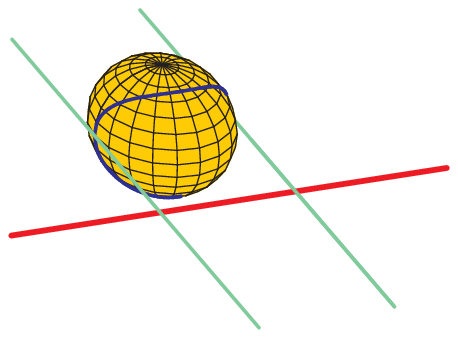}}
   \put(70,-15){(a)}
  \end{picture}
  \qquad\qquad
   \begin{picture}(150,135)(0,-10)
   \put(0,0){\epsfxsize=2in\epsffile{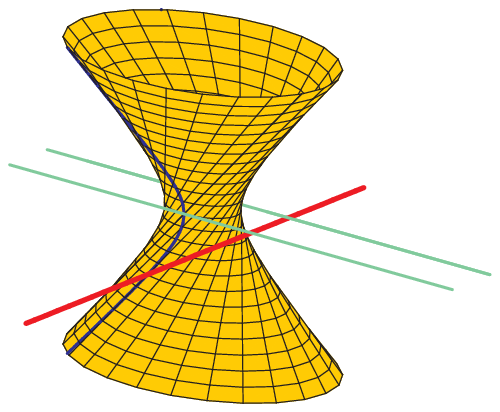}}
   \put(60,-15){(b)}
  \end{picture}
$$ 
\caption{Real quadrics not meeting the $x$-axis.\label{fig:one}}
\end{figure}
\fi
\ifpictures
\begin{figure}[htb]
$$
   \begin{picture}(150,140)(0,-15)
   \put(0,0){\epsfxsize=2in\epsffile{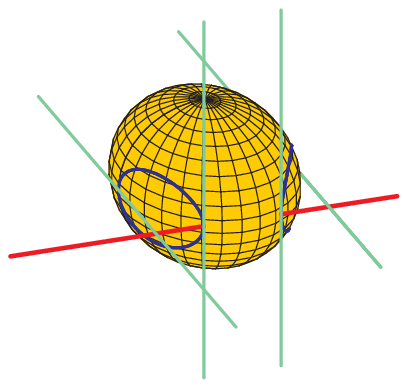}}
   \put(70,-20){(a)}
  \end{picture}
  \qquad
   \begin{picture}(150,140)(0,-15)
   \put(0,0){\epsfxsize=2in\epsffile{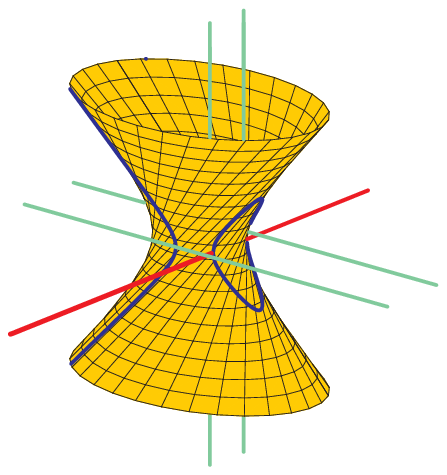}}
   \put(60,-20){(b)}
  \end{picture}
$$ 
\caption{Real quadrics meeting the $x$-axis.\label{fig:two}}
\end{figure}
\fi

\begin{rem}\label{rem:symmetry}
 For each of the spheres, there is another sphere of radius $r$ which
 leads to the same envelope, namely the one with center
 $(0,-y_0,0)^{\mathrm{T}}$.  
\end{rem}

The ramification of the $(2,2)$-curve of tangents perpendicular to the
$x$-axis is evident from Figures~\ref{fig:one} and~\ref{fig:two}.
When $x=\pm 1$, there is a single tangent line; this line has slope
$[y,z] = [1,0]$, i.e., it is a horizontal line.
When $x=\pm\sqrt{1-y_0^2}$, there is a single tangent line,
which is vertical (i.e., which has slope $[y,z]=[0,1]$). 
Figures~\ref{fig:one} and~\ref{fig:two} depict these lines in case
they are real. In Figure~\ref{fig:one} we have $|y_0| > 1$, and hence the 
vertical tangent lines are complex. 
All other values of $x$ give two lines perpendicular to the $x$-axis
and tangent to the quadric, but some have imaginary slope.
\smallskip

The difference in the number of real components of the fiber
$\varphi^{-1}(C(s))$ noted in Theorem~\ref{thm:realfamily} is evident in these
examples. 
The spheres and hyperboloid displayed together are isomorphic under the change
of coordinates $z\mapsto \sqrt{-1}\cdot z$, which interchanges the transversal
tangents of purely imaginary slope for one quadric with the real
transversal tangents of the other and corresponds to the different signs $\pm$
in~\eqref{eq:real-quad} and~\eqref{eq:spec22}. 

For the sphere of Figure~\ref{fig:one}, only 4 of the 12 families are real.
One consists of ellipsoids, including the original sphere, one of hyperboloids
of two sheets, and two of hyperboloids of one sheet.
Since a hyperboloid of two sheets can be seen as an ellipsoid meeting the plane at
infinity in a conic, we see there are two families of ellipsoids and two of
hyperboloids. 
In Figure~\ref{fig:3families}, we display one quadric from each family 
(except the family of the sphere), 
together with the original sphere, the $x$-axis, and the curve on the quadric
where the lines perpendicular to the $x$-axis are tangent to the quadric.
\ifpictures
\begin{figure}[htb]
$$
  \epsfxsize=1.8in\epsffile{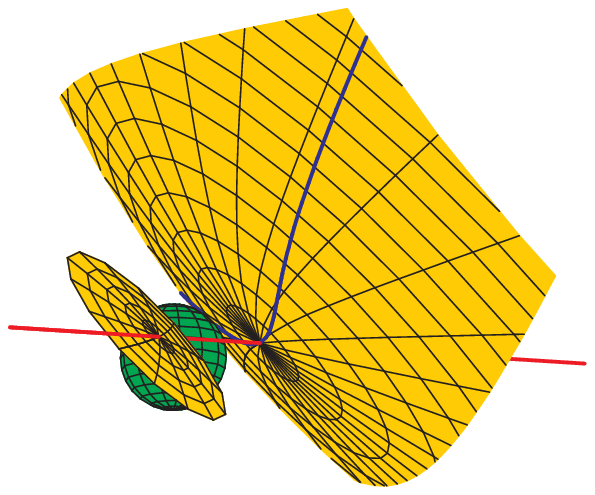}\qquad
  \epsfxsize=1.5in\epsffile{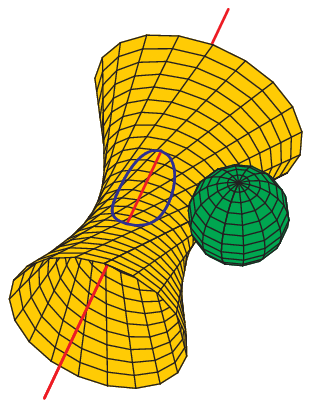}\qquad
  \epsfxsize=1.5in\epsffile{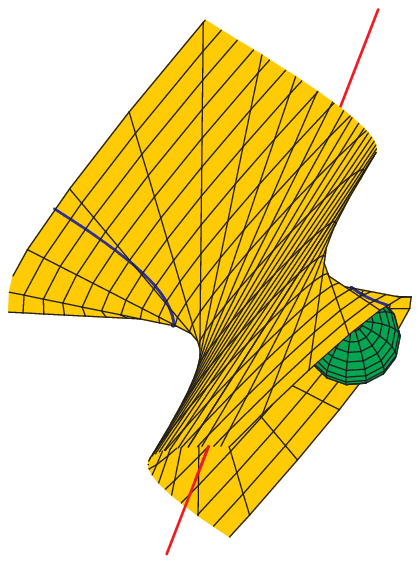}
$$
\caption{The other three families.\label{fig:3families}}
\end{figure}
\fi

Similarly, the hyperboloid of Figure~\ref{fig:one} has only 4 of its 12
families real with two families of ellipsoids and two of hyperboloids.
The sphere of Figure~\ref{fig:two} has only 4 of its 12 families real, and
all 4 contain ellipsoids.
In contrast, the hyperboloid of Figure~\ref{fig:two} has all 12 of its
families real, and they contain only hyperboloids of one sheet.

Many more pictures (in color) are found on the web page\newfootnote{{\tt
http://www.math.umass.edu/\~{}sottile/pages/2l2s/index.html}} accompanying
this article.


\section{Transversals to two lines and tangents to two
spheres}\label{se:lines2l2s} 

We solve the original question of configurations of two lines and two
spheres for which there are infinitely many real transversals to 
the two lines that are also tangent to both spheres.
While general quadrics are naturally studied in projective space $\P^3$,
spheres naturally live in (the slightly more restricted) affine space $\R^3$. 
As noted in Section~2, we treat only skew lines.
There are two cases to consider.
Either the two lines are in $\R^3$ or one lies in the plane at infinity.
We work throughout over the real numbers.

\subsection{Lines in affine space $\R^3$}  
\begin{figure}[htb]
\[
  \begin{picture}(310,90)(5,0)
   \put(0,15){\epsfxsize=1.4in\epsfbox{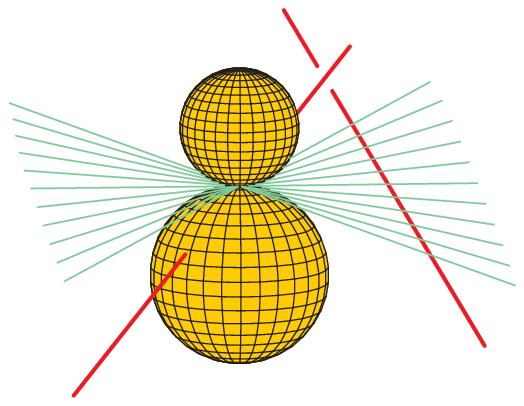}}
   \put(140,35){\epsfxsize=2.2in\epsfbox{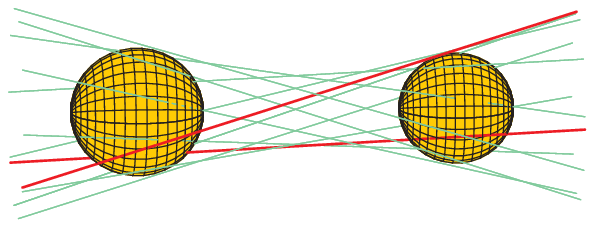}}
   \put(32,0){(1)}   \put(210,0){(2)}
  \end{picture}
\]
\caption{Examples from Theorem~\ref{th:2l2sclassification}\label{fig3}.}
\end{figure}

The complete geometric characterization of configurations where the lines lie
in $\R^3$ is stated in the following theorem and illustrated in
Figure~\ref{fig3}. 

\begin{thm}
\label{th:2l2sclassification} 
Let $S_1$ and $S_2$ be two distinct spheres 
and let $\ell_1$ and $\ell_2$ be two skew lines in $\R^3$.
There are infinitely many lines that meet $\ell_1$ and $\ell_2$ and 
are tangent to $S_1$ and $S_2$ in \emph{exactly} the following cases.
\begin{enumerate}
 \item[(1)] The spheres $S_1$ and $S_2$ are tangent to each other at
            a point $p$ which lies on one line, and the second line lies in
            the common tangent plane to the spheres at the point $p$.
            The pencil of lines through $p$ that also meet the second line 
            is exactly the set of common transversals to $\ell_1$ and $\ell_2$
            that are also tangent to $S_1$ and $S_2$. 

 \item[(2)] The lines $\ell_1$ and $\ell_2$ are each tangent to both $S_1$
            and $S_2$, and they are images of each other under a
            rotation about the line connecting the centers of $S_1$ and $S_2$.
            If we rotate $\ell_1$  about the line connecting the centers of the
            spheres, it sweeps out a hyperboloid of one sheet. 
            One of its rulings contains $\ell_1$ and $\ell_2$, and the lines
            in the other ruling are tangent to $S_1$ and $S_2$ and meet
            $\ell_1$ and $\ell_2$, except for those that are parallel to one
            of them. 
\end{enumerate}
\end{thm}

Let $\ell_1$ and $\ell_2$ be two skew lines.
The class of spheres is not invariant under the set of projective linear
transformations, but rather under the group generated by rotations,
translations, and scaling the coordinates.
Thus we can assume that
\[
  \ell_1 = \left\{ \left( \begin{matrix}
    0 \\ 0 \\ 1 
  \end{matrix} \right)
  + x
  \left( \begin{matrix}
    1 \\ \delta \\ 0 
  \end{matrix} \right) \; : \; x \in \R \right\} \,, \qquad
  \ell_2 = \left\{ \left( \begin{matrix}
    0 \\ 0 \\ -1 
  \end{matrix} \right)
  + z
  \left( \begin{matrix}
    -1 \\ -\delta \\ 0 
  \end{matrix} \right) \; : \; z \in \R \right\} \,
\]
for some $\delta\in\R \setminus \{0\}$. 
As before, there is a one-to-one correspondence between
lines meeting $\ell_1$ and $\ell_2$ and pairs $(x,z)\in\R^2$.
The transversal corresponding to a pair $(x,z)$ passes through the 
points $(x,\delta x, 1)^{\mathrm{T}}$ and 
$(z,-\delta z, -1)^{\mathrm{T}}$, and has Pl\"ucker coordinates
\[
  (z-x,-\delta(x+z),-2,-2 \delta x z,-(x+z),\delta(z-x))^{\mathrm{T}} \, .
\]
Let $S_1$ have center $(a,b,c)^{\mathrm{T}}$ and radius $r$. 
By Proposition~\ref{le:quadrictangentcond}
and~(\ref{eq:plueckertangenttosphere}), the transversals
tangent to $S_1$ are parameterized by a curve $C_1$ of degree 4
with equation
\begin{eqnarray}
0&=&4\delta^2 x^2 z^2\,+\,4\delta(b{-}a\delta)x^2z\,+\,
\bigl((b{-}a\delta)^2+(1{+}\delta^2)((1{+}c)^2-r^2)\bigr)x^2 \nonumber \\
&& -\,4\delta(b+a\delta)xz^2
\,+\,2\bigl((r^2{-}c^2)(1{-}\delta^2)+(1{-}b^2)+\delta^2(a^2{-}1)\bigr)xz 
      \label{eq:c1} \\
&&-\,4(1{+}c)(a{+}b\delta)x
\,+\,\bigl((b{+}a\delta)^2+(1{+}\delta^2)((1{-}c)^2{-}r^2)\bigr)z^2 \nonumber \\
&& +\,4(c{-}1)(a{-}b\delta)z\,+\,4(a^2+b^2-r^2)\, . \nonumber
\end{eqnarray}
This is a dehomogenized version of the bihomogeneous
equation~(\ref{eq:22}) of bidegree~$(2,2)$.
Note also that the curve $C_1$ is defined over our ground field $\R$.
The transversals to $\ell_1$ and $\ell_2$ 
tangent to $S_2$ are parameterized by a similar curve $C_2$. 
There are infinitely many lines which meet $\ell_1$ and $\ell_2$ and are
tangent to $S_1$ and $S_2$ if and only if the curves $C_1$ and $C_2$ have a
common component. 
That is, if and only if the associated polynomials share a common factor.
We first rule out the case when the curves are irreducible.

\begin{lemma}
\label{le:c1s1}
 The curve $C_1$ in~\eqref{eq:c1} determines the sphere $S_1$ uniquely. 
\end{lemma}

\begin{proof}
Given the curve~(\ref{eq:c1}), we can
rescale the equation such that the coefficient of $x^2z^2$ is
$4\delta^2$. 
 From the coefficients of $x^2z$ and $xz^2$ we can
determine $a$ and $b$, and then from the coefficients of $x^2$ and
$z^2$ we can determine $c$ and $r$. 
\end{proof}

\begin{rem}
 By Remark~\ref{rem:symmetry}, Lemma~\ref{le:c1s1} does not hold if 
 the lines are allowed to live in projective space $\P^3$. We come back
 to this in Section~\ref{se:2l2sclassificationprojective}.
\end{rem}

By Lemma~\ref{le:c1s1}, there can be infinitely many common transversals to 
$\ell_1$ and $\ell_2$ that are tangent to two spheres only if the curves 
$C_1$ and $C_2$ are reducible.
In particular, this rules out cases (1) and (2) of Lemma~\ref{le:CrClass}.
Our classification of factors of $(2,2)$-forms in Lemma~\ref{le:CrClass}
gives the following possibilities for the common irreducible factors (over
$\R$) of $C_1$ and $C_2$, up to interchanging $x$ and $z$. 
Either the factor is a cubic (the dehomogenization of a $(2,1)$-form), or it
is linear in $x$ and $z$ (the dehomogenization of a $(1,1)$-form), or it is 
linear in $x$ alone (the dehomogenization of a $(1,0)$-form).
There is the possibility that the common factor will be an irreducible (over
$\R$) 
quadratic polynomial in $x$ (coming from a $(2,0)$-form), but then this
component will have no real points, and thus contributes no common real
tangents. 

We rule out the possibility of a common cubic factor,
showing that if $C_1$ factors as $x-x_0$ and a cubic, 
then the cubic still determines $S_1$.
The vector $(-\delta,-1,\delta x_0)^{\mathrm{T}}$ is perpendicular to the
plane through $(x_0,\delta x_0,1)^{\mathrm{T}}$ and $\ell_2$,  so the center
of $S_1$ will be  
$(x_0,\delta x_0,1)^{\mathrm{T}}+\lambda (-\delta,-1,\delta x_0)^{\mathrm{T}}$
for some non-zero  $\lambda\in\R$.
Thus $r^2=\lambda^2(1+\delta^2+\delta^2 x_0^2)$. 
Substituting this into~\eqref{eq:c1} and dividing by $(x-x_0)$
we obtain the equation of the cubic:
\begin{equation}\label{cubic}
 \begin{array}{rcl}
   0&=&\delta^2 x z^2
    \,+\,\delta(\delta^2{-}1)\lambda xz 
    \,+\,(1{+}\delta^2(1{-}\lambda^2)+\delta\lambda(1{+}\delta^2)x_0)x\\
    &&+\, \delta(\lambda(1{+}\delta^2)-\delta x_0)z^2 
   \,+\,\delta(\delta^2{-}1)\lambda x_0 z
   \,+\,4\delta\lambda+(\delta^2\lambda^2-\delta^2-1)x_0\,.
   \rule{0pt}{15pt}
 \end{array}
\end{equation}
Given only this curve, we can rescale its equation so that the
coefficient of $xz^2$ is $\delta^2$, then if $\delta\neq\pm 1$, we can
uniquely determine $\lambda$, $x_0$ and therefore $S_1$, too,
from the coefficients of $xz$ and $x$.

The uniqueness is still true if $\delta=\pm1$. Assume that 
$\delta=1$. Then (\ref{cubic}) reduces to 
 \[
    xz^2\,+\,(2\lambda -x_0)z^2\,+\,(2-\lambda^2+2\lambda x_0)x
    \,+\,4\lambda+(\lambda^2-2)x_0\ =\ 0\,.
 \]
Set $\alpha:=2\lambda-x_0$, $\beta:=2-\lambda^2+2\lambda x_0$, and 
$\gamma:=4\lambda+(\lambda^2-2)x_0$. 
We can solve for $\lambda$ and $x_0$ in terms of $\alpha$ and $\beta$,
 \[
    \lambda\ =\ \frac{\alpha\pm \sqrt{\alpha^2+3\beta-6}}3, \qquad 
     x_0\ =\ \frac{-\alpha\pm 2\sqrt{\alpha^2+3\beta-6}}3\,.
 \]
(We take the same sign of the square root in both cases).
If we substitute these values into the formula for $\gamma$,
we see that the two possible values of $\gamma$ coincide if and only if 
$\alpha^2+3\beta-6=0$, in which case there is only one solution for $\lambda$
and $x_0$, so $\alpha$, $\beta$, and $\gamma$ always determine $\lambda$
and $x_0$ uniquely and hence $S_1$ uniquely. The case $\delta=-1$ is
similar.\smallskip

We now are left only with the cases when $C_1$ and $C_2$ contain a common
factor of the form $x-x_0$ or $xz+ sx + tz + u$.
Suppose the common factor is $x-x_0$.
Then any line through $p:=(x_0,\delta x_0,1)^{\mathrm{T}}$ and a point of
$\ell_2$ is tangent to $S_1$. 
This is only possible if the sphere $S_1$ is tangent to the plane through 
$p$ and $\ell_2$ at the point $p$. 
We conclude that if $C_1$ and $C_2$ have the common factor $x-x_0$, then
the spheres $S_1$ and $S_2$ are tangent to each other at
the point $p=(x_0,\delta x_0,1)^{\mathrm{T}}$ lying on $\ell_1$ and $\ell_2$
lies in the common tangent plane to the spheres at the point $p$.
This is case (1) of Theorem~\ref{th:2l2sclassification}.

Suppose now that $C_1$ and $C_2$ have a common irreducible
factor $xz+ sx + tz + u$.  
We can solve the equation $xz+ sx + tz + u=0$
uniquely for $z$ in terms of $x$ for general values of $x$, or for $x$ in
terms of $z$ for general values of $z$, this
gives rise to an isomorphism $\phi$ between the
projectivizations of $\ell_1$ and $\ell_2$. 
The lines connecting $q$ and $\phi(q)$
as $q$ runs through the points of $\ell_1$ sweep out a hyperboloid of
one sheet. 
The lines $\ell_1$ and $\ell_2$ are contained in one ruling, and the lines
meeting both of them and tangent to $S_1$ are the lines in the other
ruling. 

\begin{lemma}
 Let $H\subset\R^3$ be a hyperboloid of one sheet. If all lines
 in one of its rulings are tangent to a sphere $S$, then $H$ is a
 hyperboloid of revolution, the center of the sphere $S$ is on the axis of
 rotation and $S$ is tangent to $H$.
\end{lemma}

\begin{proof}
We can choose Cartesian co-ordinates such that $H$ has equation 
$x^2/\alpha^2+y^2/\beta^2-z^2/\gamma^2=1$ for some positive real
numbers $\alpha$, $\beta$, $\gamma$.
Let the sphere have center
$(A,B,C)^{\mathrm{T}}$ and radius $R$. The set of points of the form 
$(\alpha,\lambda \beta,\lambda \gamma)^{\mathrm{T}}$, 
$(-\alpha,-\lambda \beta,\lambda \gamma)^{\mathrm{T}}$,
$(-\lambda \alpha, \beta,\lambda \gamma)^{\mathrm{T}}$ and 
$(\lambda \alpha,- \beta,\lambda \gamma)^{\mathrm{T}}$ as
$\lambda$ runs through $\R$ form four lines in one of the rulings.
Since the two rulings are symmetric, we only need to deal with one of
them. 

The sphere $S$ is tangent to a line if and only if the distance of the
line from the center of $S$ is $R$.
The condition that $(A,B,C)^{\mathrm{T}}$ must be at the same distance from
the first two lines gives the equation
\[
  \alpha(\beta^2+\gamma^2)A\,+\,\beta\gamma BC\ =\ 0\,,
\]
the equality of distances from the othe two lines gives
\[
  \beta(\alpha^2+\gamma^2)B\,-\,\alpha\gamma AC\ =\ 0\,.
\]

Since $\alpha,\beta,\gamma>0$, 
the common solutions of these equations have $A=B=0$. 
Using this
information, the equality of the distances from the first and third
lines gives $\alpha=\beta$, or $C=\pm\sqrt{(a^2+c^2)(b^2+c^2)/c}$. 
To eliminate this second possibility, consider two more lines
in the same ruling, the points of the form
\[
  \Bigl(\alpha\frac{1-\lambda}{\sqrt 2},\ \beta\frac{1+\lambda}{\sqrt 2},\ 
                     \lambda\gamma\Bigr)
   \quad\textrm{ and }\quad
  \Bigl(\alpha\frac{1+\lambda}{\sqrt 2},\ \beta\frac{-1+\lambda}{\sqrt 2},\ 
    \lambda\gamma\Bigr)
\]
as $\lambda$ runs through $\R$. The equality of distances from these
two lines together with $A=B=0$ gives $\alpha=\beta$ or $C=0$.

Therefore the only case when $(A,B,C)^{\mathrm{T}}$ can be  at the same
distance from all lines in one ruling of $H$ is when $\alpha=\beta$, i.e., $H$
is a hyperboloid of revolution about the $z$-axis, and $(A,B,C)^{\mathrm{T}}$
lies on the $z$-axis. 
In this case, it is obvious that $(A,B,C)^{\mathrm{T}}$ is at the same
distance from all the lines contained in $H$, and these lines are
tangent to $S$ if and only if $S$ is tangent to $H$.
\end{proof}

By this lemma, the hyperboloid swept out by the lines
meeting $\ell_1$ and $\ell_2$ and tangent to $S_1$
is a hyperboloid of revolution with the center of $S_1$ 
on the axis of rotation.
Furthermore, $\ell_1$ and $\ell_2$ are lines in one the rulings
of the hyperboloid, therefore they are images of each other under
suitable rotation about the axis, the images of $\ell_1$ sweep out the
whole hyperboloid, and $\ell_1$, $\ell_2$ are both tangent to $S_1$. 
Applying the lemma to $S_2$ shows that the center of $S_2$ is also on
the axis of rotation and $\ell_1$, $\ell_2$ are both tangent to $S_2$. 
We cannot have $S_1$ and $S_2$ concentric, therefore the axis of rotation
is the line through their centers.
This is exactly case (2) of Theorem~\ref{th:2l2sclassification}, and we
have completed its proof.

\subsection{Lines in projective space. \label{se:2l2sclassificationprojective}}
We give the complete geometric characterization of configurations
in real projective space  
where the line $\ell_2$ lies in the plane at infinity.

\begin{thm}
\label{th:2l2sclass} 
Let $S_1$ and $S_2$ be two distinct spheres 
and let $\ell_1$ lie in  $\R^3$ with 
$\ell_2$ a line at infinity skew to $\ell_1$.
There are infinitely many lines that meet $\ell_1$ and $\ell_2$ and 
are tangent to $S_1$ and $S_2$ in \emph{exactly} the following cases.
\begin{enumerate}
 \item[(1)] The spheres $S_1$ and $S_2$ are tangent to each other at
            a point $p$ which lies on $\ell_1$, and $\ell_2$ is the line at
            infinity in the common tangent plane to the spheres at the point
            $p$. 
            The pencil of lines through $p$ that lie in this tangent plane
            are exactly the common transversals to $\ell_1$ and $\ell_2$ that
            are also tangent to $S_1$ and $S_2$. 
 \item[(2)] Any line meeting $\ell_1$ and $\ell_2$ is perpendicular to
            $\ell_1$ and $S_1$ and $S_2$ are related to each other by
            multiplication by $-1$ in the directions perpendicular to
            $\ell_1$.
            Thus we are in exactly the situation of Remark~\ref{rem:symmetry}
            of Section~\ref{se:examples} as shown in Figures~\ref{fig:one}(a)
            and~\ref{fig:two}(a).
\end{enumerate}
\end{thm}

\begin{proof}
Let $\Pi$ be any plane passing through a point of $\ell_1$ and containing
$\ell_2$.
Then common transversals to $\ell_1$ and $\ell_2$ are lines meeting $\ell_1$
that are parallel to $\Pi$.
Choose a Cartesian coordinate system in $\R^3$ such that $\ell_1$ is the
$x$-axis.
Suppose that $S_1$ has center $(a,b,c)^{\mathrm{T}}$ and radius $r$. 
Let $u=(u_1,u_2,0)^{\mathrm{T}}$ and  $v=(v_1,0,v_3)^{\mathrm{T}}$ be vectors
with $u_2\neq 0$ and $v_3\neq 0$ parallel to $\Pi$. 
Such vectors exist as $\ell_1$ and $\ell_2$ are skew.
A common transversal to $\ell_1$ and $\ell_2$ is determined
by the intersection point $(x,0,0)^{\mathrm{T}}$ with $\ell_1$ and a direction
vector corresponding to the intersection point with $\ell_2$, which can be
written as $u+zv$ for some $z\in \R$, unless it is parallel to $v$. 
Since $S_1$ has at most two tangent lines which meet $\ell_1$ that are
parallel to $v$, so by omitting these we are not 
losing an infinite family of common transversals/tangents.

The transversals that are tangent to $S_1$ are parametrized by a curve
$C_1$ in the $xz$-plane with equation 
\begin{eqnarray}
 0&=&v_3^2 x^2 z^2+u_2^2 x^2+ 2 v_3(c v_1 - a v_3) x z^2+
 2(b u_2 v_1 +c u_1 v_3)xz \nonumber \\
 &&+2 u_2(b u_1-\!au_2)x+((b^2+c^2-\! r^2)v_1^2-\! 
 2acv_1v_3+(a^2+b^2-\!r^2)v_3^2)z^2 \label{c1'} \\
 &&+2((b^2+c^2-r^2)u_1v_1-acu_1v_3-bu_2(av_1+c v_3))z \nonumber \\
 &&+((b^2+c^2- r^2)u_1^2-2abu_1u_2+(a^2+c^2-r^2)u_2^2) \nonumber
\end{eqnarray}
The transversals tangent to $S_2$ are parametrized by a similar curve $C_2$. 
There are infinitely many lines that meet $\ell_1$ and $\ell_2$ and are
tangent to $S_1$ and $S_2$ if and only if  $C_1$ and $C_2$ have a common
non-empty real component. 

It is easy to see from the coefficients of $xz^2$, $xz$ and $x$ and
the constant term that if $u_1\neq 0$ or $v_1\neq 0$, then $C_1$  
determines $a$, $b$, $c$ and $r^2$ and therefore $S_1$ uniquely, so if 
$C_1$ is irreducible and $u_1\neq 0$ or $v_1\neq0$, then there cannot 
be infinitely many common transversals that are tangent to $S_1$ and
$S_2$.

Assume now that  $u_1=v_1=0$, this is equivalent to the plane $\Pi$
being perpendicular to $\ell_1$. 
 From the coefficient of $x$ we can determine $a$, and then from the
coefficients of $z^2$, $z$, and the constant term we can calculate the
quantities $\alpha=c^2-r^2$,  $\beta=bc$, and $\gamma=b^2-r^2$. 
The equation $(\alpha+r^2)(\gamma+r^2)-\beta^2=0$ is a quadratic equation for
$r^2$ with solutions 
\[
  r^2\ =\ \frac{1}{2}
   \bigl(-\alpha-\gamma\pm\sqrt{(\alpha-\gamma)^2+4\beta^2}\bigr)\,.
\]
Only the larger root is feasible, even when both are positive, since both
$\alpha+r^2=c^2$  and $\gamma+r^2=b^2$ must be non-negative. 
Hence $r^2$, and thus $b^2$ and $c^2$ are uniquely determined. 
The values of $b^2$, $bc$, and $c^2$ determine two possible pairs 
$(b,c)$ which are negatives of each other. 
This is exactly case (2) of the theorem.
In fact, this case is illustrated by Figures~\ref{fig:one}(a)
and~\ref{fig:two}(b).

Let us now consider the cases when $C_1$ is reducible. 
As in the proof of Theorem~\ref{th:2l2sclassification},
we need only consider cubics and factors of the form
$xz+sx+tz+u$, $x-x_0$, and $z-z_0$. 

Assume that $C_1$ has a component with equation $xz+sx+tz+u$. As
described in the proof of Theorem~\ref{th:2l2sclassification}, this
establishes an isomorphism between the projectivizations of $\ell_1$
and $\ell_2$. The lines connecting the corresponding points of the
projectivizations of $\ell_1$ and $\ell_2$ sweep out a hyperbolic
paraboloid. However, the lines in one ruling of the hyperbolic
paraboloid cannot all be tangent to a sphere, therefore this case
cannot occur.

Likewise, the factor $z-z_0$ cannot appear, since it would mean that 
all the lines through a point of $\ell_1$ parallel to a certain
direction are tangent to $S_1$, which is clearly impossible.

Consider the case where the equation of $C_1$ has a factor of $x-x_0$.
As we saw in the proof of Theorem~\ref{th:2l2sclassification},
$\ell_1$ meets the sphere $S_1$ at the point $p:=(x_0,0,0)^{\mathrm{T}}$, 
and $\ell_2$ lies in the
tangent plane to $S_1$ at $p$, and so this tangent plane is parallel
to $\Pi$.

If $x-x_0$ is a factor of $C_2$, too, then $C_2$ passes through $p$
and its tangent plane there is also parallel to $\Pi$, so we have case
(1) of the theorem.

To finish the proof we investigate what happens 
if the common component of $C_1$ and $C_2$ is the cubic 
obtained from $C_1$ after removing the line 
$x-x_0=0$.

The center of $S_1$ has
coordinates $(x_0+\mu u_2 v_3,-\mu u_1 v_3,-\mu u_2 v_1)^{\mathrm{T}}$ for
some $\mu\in\R$, since $S_1$ passes through $(x_0,0,0)^{\mathrm{T}}$ and its
tangent plane there is parallel to $\Pi$, and we have 
$r^2=\mu^2 (u_1^2 v_3^2 + u_2^2 v_1^2+u_2^2 v_3^2)$. Substituting
this into (\ref{c1'}) we obtain the equation of the remaining cubic,
\begin{eqnarray*}
&&v_3^2 x z^2+u_2^2 x -v_3(x_0 v_3+2\mu u_2(v_1^2+v_3^2))z^2\\
&&-4\mu u_1 u_2v_1v_3 z-u_2(x_0u_2+2\mu v_3(u_1^2+u_2^2))=0.
\end{eqnarray*}
If $u_1\neq 0$ or $v_1\neq 0$ then from the coefficients of this curve
we can determine $x_0$ and $\mu$, hence $S_1$ uniquely, so $C_1$ and
$C_2$ cannot have a common cubic component.
If $u_1=v_1=0$ then the above equation factorizes as
$$(x- (2 \mu u_2 v_3  + x_0))(v_3^2 z^2+u_2^2)=0,$$
so if $C_2$ contains the curve defined by this equation, then the line
$x- (2 \mu u_2 v_3  + x_0)=0$ is a common component of both $C_1$ and
$C_2$, which is a case we have already dealt with.
\end{proof}

\appendix
\section{Calculations from Section \ref{se:infinitetangents}}

We describe the computation of Section~\ref{se:infinitetangents} 
in much more detail, giving a 
commentary on the {\sc Singular} file that accomplishes the computation and
displaying its output. 
The input and output are displayed in {\tt typewriter} font on separate lines
and the output begins with the {\sc Singular} comment characters 
({\tt //}).

The library {\tt primdec.lib} contains the
function {\tt sat} for saturating ideals, and the option {\tt redSB}
forces {\sc Singular} to work with reduced Gr\"obner (standard) bases.
\begin{verbatim}
  LIB "primdec.lib";
  option(redSB);
\end{verbatim}
We initialize our ring.
\begin{verbatim}
  ring R = 0, (s,t,a,b,c,d,e,f,g,h,k,l), (dp(2), dp(10));
\end{verbatim}
The underlying coefficient field has characteristic 0 
(so it is $\mathbb{Q}$) and variables
$s,t,a,\ldots,k,l$, with a product term order 
chosen to simplify our analysis of the projection 
to $\mathbb{C}^2$, the space of parameters. 

We consider the ideal generated by~\eqref{eq:quad-gens}
\begin{verbatim}
  ideal I = el-g^2, ek-gf, ak-dc, ah-c^2;
\end{verbatim}
and by 
the $2\times 2$ minors of the coefficient matrix~\eqref{eq:coeff-matrix}.
\begin{verbatim}
  matrix M[2][5] =    s   ,    1-s   ,       -2     ,     1-t  ,   t   , 
                    al-d^2, 2*(bl-dg), 2*(2bk-cg-df), 2*(bh-cf), eh-f^2; 
  I = I + minor(M,2);
\end{verbatim}
We check the dimension and degree (multiplicity) of the variety
$\calV(I)$, first computing a Gr\"obner basis for $I$.
\begin{verbatim}
  I = std(I);  dim(I), mult(I);
  //   6  8
\end{verbatim}
{\sc Singular} gives the dimension of $\calV(I)$ in affine space
$\mathbb{C}^{12}$.
Since $I$ is homogeneous in the variables $a,b,\ldots,h,k,l$, 
we consider $\calV(I)$ to be a subvariety of 
$\mathbb{P}^9\times\mathbb{C}^2$.
Its dimension is one less than that of the corresponding affine variety.
Thus $\calV(I)$ has dimension 5 and degree 8.

In Section~\ref{se:infinitetangents}, we identified three
spurious components of $\calV(I)$ which we remove.
The first and largest is the ideal of rank 1 quadrics, given by the 
$2\times 2$-minors of the $4\times 4$-symmetric matrix~\eqref{eq:symmat}.
\begin{verbatim}
  matrix Q[4][4] =  a , b , c , d , 
                    b , e , f , g , 
                    c , f , h , k , 
                    d , g , k , l ;
  ideal E1 = std(minor(Q,2));
\end{verbatim}
We remove this spurious component, computing the quotient ideal 
$(I\;:\, E_1)$.
\begin{verbatim}
  I = std(quotient(I,E1));  dim(I), mult(I);
  //   5  20
\end{verbatim}
The other two spurious components describe rank 2 quadrics which are unions
of two planes with intersection line $\ell_1$ or $\ell_2$.
\begin{verbatim}
  ideal E2 = g, f, e, d, c, b, a;  // intersection line l1
  ideal E3 = l, k, h, g, f, d, c;  // intersection line l2
\end{verbatim}
The corresponding components are not reduced; rather than take ideal
quotients, we saturate the ideal $I$ with respect to these spurious ideals.
The {\sc Singular} command {\tt sat} for saturation returns a pair whose first
component is a Gr\"obner basis of the saturation and the second is the
saturation exponent. Here, both saturations have exponent 4.
We saturate $I$ with respect to $E_2$,
\begin{verbatim}
  I = sat(I,E2)[1];  dim(I), mult(I);
  //   5  10
\end{verbatim}
and then with respect to $E_3$.
\begin{verbatim}
  ideal J = sat(I,E3)[1];  dim(J), mult(J);
  //   4  120
\end{verbatim}
Thus we now have a variety $\calV(J)$ of dimension 3 in
$\mathbb{P}^9\times\mathbb{C}^2$. 
We check that it projects onto the space $\mathbb{C}^2$ of parameters by
eliminating the variables $a,b,\ldots,h,k,l$ from $J$.
\begin{verbatim}
  eliminate(J, abcdefghkl);
  //   _[1]=0
\end{verbatim}
Since we obtain the zero ideal, the image of $\calV(J)$ is 
Zariski dense in
$\mathbb{C}^2$ \cite[Chapter~4, \S 4]{clo-b96}.
However, the projection 
$\mathbb{P}^9\times\mathbb{C}^2\twoheadrightarrow\mathbb{C}^2$ is a closed
map, so the image of $\calV(J)$ is $\mathbb{C}^2$.
Thus, for
every smooth $(2,2)$-curve $C$ defined by~(\ref{eq:gen22}), there is a 
quadric whose transversal tangents are described by the curve $C$.

We now apply the factorization Gr\"obner basis algorithm {\tt facstd} to $J$.
The second argument of {\tt facstd} is the list of non-zero constraints
which are given by Theorem~\ref{thm:gen22}.
\begin{verbatim}
  ideal L = s, t, t-1, s-1, s-t; 
  list F = facstd(J,L);
\end{verbatim}
{\sc Singular} computes seven factors
\begin{verbatim}
  size(F);
  //   7
\end{verbatim}
Since $J$ and the seven factors $L_1, \ldots, L_7$ are radical ideals,
this factorization can be verified by checking that
that the following ideals $\mathit{V}_1$ and $\mathit{V}_2$ coincide.
(This part of the calculation is archived on the web 
page\newfootnote{{\tt http:/$\!$/www.math.umass.edu/\~{}sottile/pages/2l2s.html}}.)
\begin{verbatim}
  int i;
  ideal FF = 1;
  for (i = 1; i <= 7; i++) { FF = intersect(FF,F[i]); }
  ideal V1, V2;
  V1 = std(sat(sat(sat(sat(sat(FF,t)[1],s)[1],t-1)[1],s-1)[1],s-t)[1]);
  V2 = std(sat(sat(sat(sat(sat(J ,t)[1],s)[1],t-1)[1],s-1)[1],s-t)[1]);
\end{verbatim}

Note, in particular, that for any given explicit values of $s,t$ satisfying
the nonzero conditions, the parametric factorization (in $s,t$) 
produced by {\tt facstd} can be specialized to an explicit factorization.

We examine the ideals in the list $F$, working over the ring with parameters.
\begin{verbatim}
  ring S = (0,s,t), (a,b,c,d,e,f,g,h,k,l), lp; short = 0;
\end{verbatim}
First, the ideal $J$ has dimension 1 and degree 24 over this ring,
as claimed. 
\begin{verbatim}
  ideal JS = std(imap(R,J));   dim(JS), mult(JS);
  //   2  24
\end{verbatim}
The first ideal in the list $L$ has dimension 0.
\begin{verbatim}
  setring R; FR = F[1]; setring S; 
  FS = std(imap(R,FR));  dim(FS), mult(FS);
  //   1  4
\end{verbatim}
This ideal is a spurious component from the factorization. It is
contained in the spurious ideal $E_2$.
\begin{verbatim}
  FS[5], FS[6], FS[7], FS[8], FS[9], FS[10], FS[11];
  //   g f e d c b a
\end{verbatim}
The other six components each have dimension 1 and degree 4, and each contains
a homogeneous quadratic polynomial in the variables $x$ and $y$.
\begin{verbatim}
  for (i = 2; i <= 7; i++) {
    setring R;  FR = F[i];  setring S; 
    FS = std(imap(R,FR));  dim(FS), mult(FS);
    FS[1];
    print("--------------------------------");
  }
  //   2  4
  //   (-s^2+2*s-1)*k^2+(2*s-2)*k*l+(s*t-1)*l^2
  //   --------------------------------
  //   2 4
  //   (s-1)*k^2-2*k*l-l^2
  //   --------------------------------
  //   2 4
  //   (s^2-2*s+1)*k^2+(-2*s+2)*k*l+(-t+1)*l^2
  //   --------------------------------
  //   2 4
  //   (s^2-2*s+1)*k^2+(-2*s+2)*k*l+(-t+1)*l^2
  //   --------------------------------
  //   2 4
  //   (s-1)*k^2-2*k*l-l^2
  //   --------------------------------
  //   2 4
  //   (-s^2+2*s-1)*k^2+(2*s-2)*k*l+(s*t-1)*l^2
  //   --------------------------------
\end{verbatim}

The whole computation takes 7 minutes CPU time on an 800 Mhz Pentium III
processor, and 3 minutes of that time are spent on the {\tt facstd} operation.

Each of these homogeneous quadratic polynomials factors over 
$\mathbb{Q}(\sqrt{s},\sqrt{t})$, and induces a factorization of
the corresponding ideal. We describe this factorization---which is
carried out by hand---in detail
for the second component $F_2$. We start from the Gr\"obner basis
of the ideal $F_2$ computed in the program above,
\begin{equation}
  \label{eq:generators}
  \begin{array}{c}
    (s-1)k^2 - 2kl - l^2, \
    (s-1)h + (2t-2)k + (t-1)l, \
    fl - gk, \\ [0.3ex]
    el-g^2, \
    d+f+g, \
    c, \
    2b+e, \ 
    a, \\ [0.3ex]
    (s-1)fk - 2gk - gl, \
    (s-1)f^2 - 2fg - g^2, \
    ek-fg \, .
  \end{array}
\end{equation}
Over $\mathbb{Q}(\sqrt{s},\sqrt{t})$, the first polynomial factors into
\[
  \left( (\sqrt{s}+1)k + l \right) \left( (\sqrt{s}-1)k - l \right) \, .
\]
We consider the first factor; the second one can be treated similarly.
Substituting $l = - (\sqrt{s}+1)k$ into the generator $fl-gk$,
that one factors into
\[
  -k  \left( (\sqrt{s} + 1) f + g \right) \, .
\]
Since any zero of $F_2$ with $k=0$ would imply $a=c=d=f=g=h=k=l=0$ 
and thus be contained in $\mathcal{V}(E_3)$, we can
divide by $k$ and obtain a linear polynomial.
Altogether, the first two rows of~\eqref{eq:generators} become 
a set of seven independent linear 
polynomials and one quadratic polynomial $el-g^2$.
For any pair $(s,t)$ satisfying~\eqref{eq:discrimin} they define a plane conic.
We leave it to the reader to verify that the three polynomials in the 
third row are contained in the ideal generated by the first two rows.

In order to show that for none of the parameters $s$, $t$ 
satisfying~\eqref{eq:discrimin} the map $\varphi$ vanishes identically 
on this conic, consider the following point $p$ on it: 
\[
\begin{array}{c}
  (0, \
   -(\sqrt{s}+1)(s-1), \ 
   0, \
   -2\sqrt{s}(s-1), \
   2(\sqrt{s}+1)(s-1), \ 
   -2(s-1), \\ [0.3ex]
   2(\sqrt{s}+1)(s-1), \
   4(t - 1)-2 (t-1)(\sqrt{s}+1), \ 
   -2(s-1), \
   2(\sqrt{s}+1)(s-1))^2 \, .
\end{array}
\]
The coefficient of $w^2 z^2$ in $\varphi(C)$ is
\[
   -4 s (\sqrt{s}-1)^2 (\sqrt{s}+1)^2 \, ,
\]
so $\varphi(C)$ does not not vanish identically.

In order to show that for all parameters $s$, $t$ 
satisfying~\eqref{eq:discrimin} the 12~conics are distinct,
consider the quadratic polynomials in $k$ and $l$
in the {\sc Singular} output above. In the factorization
over $\mathbb{Q}(\sqrt{s},\sqrt{t})$, the ideal of each of the 12 conics
contains a generator which is linear in $k$ and $l$ and independent of 
$a,\ldots, h$. 
To show the distinctness of two conics, we distinguish two cases.

If these linear homogeneous polynomials are distinct 
(over $\mathbb{Q}(s,t)$), then it can be checked that for every given 
pair $(s,t)$ they define subspaces whose restrictions to
$(k,l) \neq (0,0)$ are disjoint.

In case that the linear homogeneous polynomials coincide then it can be
explicitly checked that both conics are distinct.
For example, both $F_2$ and $F_5$ contain the factor 
$(\sqrt{s}+1)k + l$ in the first polynomial. 
As seen above, the corresponding conic of 
$F_2$ is contained in the subspace $a=c=0$. Similarly, the corresponding
conic of $F_5$ is contained in $e=g=0$. Assuming that the two conics
are equal for some pair $(s,t)$, the equations of the ideals can be 
used to show further $a=b=c=\dotsb=h=0$.
However, due to the saturation with the excess component $E_2$ this is
not possible, and hence the two conics are distinct.

The same calculations for the other components are archived at the web page of
this paper, 
{\tt http:/$\!$/www.math.umass.edu/\~{}sottile/pages/2l2s.html}.

\section*{Acknowledgements}
We thank Dan Grayson and Mike Stillman, who, as editors for the book
``Computations in Algebraic Geometry with Macaulay $2$''
introduced the typography used in the Appendix.

\providecommand{\bysame}{\leavevmode\hbox to3em{\hrulefill}\thinspace}
\providecommand{\MR}{\relax\ifhmode\unskip\space\fi MR }
\providecommand{\MRhref}[2]{%
  \href{http://www.ams.org/mathscinet-getitem?mr=#1}{#2}
}
\providecommand{\href}[2]{#2}

\end{document}

\bibliographystyle{amsplain}
\bibliography{bibl}
 
\end{document}